\documentclass[11pt, a4paper, oneside,leqno]{article}
\usepackage{amssymb,amsmath, amsfonts, amsthm,dsfont,enumerate,graphicx}
\newtheorem{mylem}{Lemma}[section]

\newtheorem{remark}{Remark}[section]

\newcommand*{\approxident}{%
  \mathrel{\vcenter{\offinterlineskip
  \hbox{$\sim$}\vskip-.35ex\hbox{$\sim$}\vskip-.35ex\hbox{$\sim$}}}}

\title{A multi-term basis criterion for families of dilated periodic functions}

\author{Lyonell Boulton$^1$ \quad \text{and} \quad Houry Melkonian$^2$ \\ \ \\ {\small \it Department of Mathematics and} \\
{\small \it Maxwell Institute for Mathematical Sciences} \\ {\small \it Heriot-Watt University, Edinburgh EH14 4AS, UK.}}
\date{October 2017}

\begin{document}
\maketitle
\footnotetext[1]{Email address: \texttt{L.Boulton@hw.ac.uk}}
\footnotetext[2]{Email address: \texttt{hm189@hw.ac.uk}}

\begin{abstract}
In this paper we formulate a concrete method for determining whether a system of dilated periodic functions forms a Riesz basis in $L^2(0,1)$. This method relies on a general framework developed by Hedenmalm, Lindqvist and Seip about 20 years ago, which turns the basis question into one about the localisation of the zeros and poles of a corresponding analytic multiplier. Our results improve upon various criteria formulated previously, which give sufficient conditions for invertibility of the multiplier in terms of sharp estimates on the Fourier coefficients. Our focus is on the concrete verification of the hypotheses by means of analytical or accurate numerical approximations. We then examine the basis question for profiles in a neighbourhood of a non-basis family generated by periodic jump functions. For one of these profiles, the $p$-sine functions, we determine a threshold for positive answer to the basis question which improves upon those found recently. 
\end{abstract}



\newpage

\section{Introduction}
Let $f:\mathbb{R}\longrightarrow \mathbb{C}$ be a 2-periodic function such that $f\in L^2(0,1)\equiv L^2$. Consider dilations $f_n(x)=f(nx)$ and set $\mathfrak{E}_f=\{f_n\}_{n=1}^\infty$. Let $\{g_n\}_{n=1}^{\infty}\subset L^2$ be another sequence. If there exists a linear homeomorphism $T:L^2\longrightarrow L^2$ such that $Tg_n=f_n$, then $\mathfrak{E}_{f}$ and $\{g_n\}_{n=1}^{\infty}$ are said to be fully equivalent. We write $\mathfrak{E}_{f}\approxident \{g_n\}_{n=1}^{\infty}$. This relation is an equivalence relation between sequences in $L^2$, which preserves the different notions of bases (for Hilbert and also Banach spaces \cite[\S I.8]{Singer1970}). In particular, let $\operatorname{s}(x)=\sin(\pi x)$. The Fourier family $\mathfrak{E}_{\operatorname{s}}$ is an orthonormal basis of $L^2$. If $\mathfrak{E}_{f}\approxident \mathfrak{E}_{\operatorname{s}}$, then $\mathfrak{E}_f$ is a Riesz basis of $L^2$. The main purpose of this paper is to examine a general criterion for determining whether $\mathfrak{E}_f$ is fully equivalent to the Fourier basis. Our emphasis is on the concrete verifiability of the hypothesis for $f$ given explicitly, rather than on the formulation of an abstract principle.

Let the (sine) Fourier coefficients of $f$ be
\[
     \hat{f}(j)=2\int_0^1 f(x) \sin(j \pi x) \mathrm{d}x \qquad \qquad j\in \mathbb{N}.
\]
Let 
\[
     m_{f}(z)=\sum_{j=1}^\infty \frac{\hat{f}(j)}{j^{z}} 
\]
be the associated  Dirichlet series (multiplier) originally defined for 
$\Re(z)>\frac12$. By virtue of a framework developed by Hedenmalm, Lindqvist and Seip about 20 years ago \cite{HLS1997,HLS1999}, the basis question for $\mathfrak{E}_f$ can be recast in terms of questions on the localisation of the zeros and poles of $m_f$. Indeed, according to  \cite[Theorem~3.1]{HLS1997}, $\mathfrak{E}_f\approxident \mathfrak{E}_{\operatorname{s}}$
 if and only if $m_f(z)$ extends to an analytic function  which is bounded and away from zero for $\Re(z)>0$. That is
 \[
      \sup_{\Re(z)>0}|m_{f}(z)|<\infty \quad \text{and} \quad
     \inf_{\Re(z)>0}|m_{f}(z)|>0.
 \] 

Let
\[
     \operatorname{J}(x)= \operatorname{sign}(\sin(\pi x))
\]
We know that $\mathfrak{E}_{\operatorname{J}}$ is not fully equivalent to a basis of $L^2$ because it is not total. However, there exist functions with profile arbitrarily close to that of 
$\operatorname{J}(x)$ such that the corresponding dilations form a Riesz basis. In order to see this, we modify slightly an example found in \cite[p.28]{HLS1997}. Below and elsewhere, $\zeta(z)$ denotes the Riemann zeta function. For $\varepsilon\geq 0$, let
\[
    a_{\varepsilon}(j)=\begin{cases} 0 & \forall j\equiv_2 0 \\
\frac{4}{\pi j^{1+\varepsilon}} & \forall j\equiv_2 1 \end{cases}
\qquad
\text{and}
\qquad
    \operatorname{J}_{\varepsilon}(x)=\sum_{j=1}^\infty a_{\varepsilon}(j)\sin(j\pi x).
\]
Then $\operatorname{J}_0(x)=\operatorname{J}(x)$.
Since
\[
     |a_{0}(j)-a_{\varepsilon}(j)|=\frac{4}{\pi j}
     \left(1-\frac{1}{j^{\varepsilon}}\right)
\]
and the right side of this is increasing in $\varepsilon$, the monotone convergence theorem yields
\[
     \lim_{\varepsilon\to 0}\sum_{j=1}^\infty|a_{0}(j)-a_{\varepsilon}(j)|^2= 0 .
\]
Hence
 $\operatorname{J}_{\varepsilon}\to 
\operatorname{J}$ in $L^2$. 
Now, the multiplier associated to $\operatorname{J}_{\varepsilon}(x)$ is
\[
      m_{\operatorname{J}_{\varepsilon}}(z)= 
   \frac{4}{\pi}\left[1-\frac{1}{2^{1+z+\varepsilon}}\right]\zeta(1+z+\varepsilon)
\]
which has all its zeros in $\Re(z)<-\varepsilon$ and a single pole at $z=-\varepsilon$. Hence
\[
    \mathfrak{E}_{\operatorname{J}_{\varepsilon}}\approxident \mathfrak{E}_{\operatorname{s}}
\]
for all $\varepsilon>0$. Therefore there are indeed functions arbitrarily close to $\operatorname{J}(x)$ (in $L^2$ norm), with dilations a Riesz basis of $L^2$. 

Despite of the above fundamental criterion and this example, in practice it can be very difficult to determine whether $m_{f}(z)$ is bounded and away from zero, even for simple profile functions $f$ (think of the Riemann hypothesis and see \S\ref{section3}-\ref{section5} below). In \S\ref{section2} we identify sufficient conditions for the multiplier to be invertible, in terms of  $|\hat{f}(j)|$. The actual statement and proof are elementary, but our emphasis here is on the computability of all the quantities involved. Our statement (Lemma~\ref{cri3}) extends those formulated in \cite[\S4]{BBCDG2006} and \cite[\S4 and \S7]{BL2014}, which have proven to be useful for determining bases properties for the $p$-trigonometric functions.  

 By ``computable'' we mean that the hypotheses are not just abstract or given ``in principle'', but rather they can be verified for concrete profile functions $f$ by either analytical or accurate numerical means in a finite (small) number of steps. In the subsequent sections \S\ref{section3}-\ref{section5}, we derive full equivalence to the Fourier basis for three profiles in a regime very close to that of $\operatorname{J}(x)$.

The profile discussed in \S\ref{section5} is the $p$-sine functions. The full equivalence question for these functions has received significant attention in recent years \cite{BBCDG2006,BushellEdmunds2012,BL2014}, as they play a fundamental role in Approximation Theory, in the particular context of Sobolev embeddings, \cite{LangEdmunds2011}.

Below we report on various analytical and numerical thresholds.  When we display numerical quantities, these are accurate to the 6th significant figure shown and the last digit has been rounded.
We have computed all these numerical quantities with an accuracy of 12 digits or more. 

We include various results involving the Fourier coefficients of the $p$-sine functions in an Appendix. These can be regarded as independent from the rest of the text. 


\section{The multi-term criterion} \label{section2}

Let $\mathbb{P}(\mathbb{N})\subset\mathbb{N}$ be the set of all prime numbers not including 1. Let $\mathcal{F}\subset \mathbb{N}$ be a finite set such that $1\in \mathcal{F}$. Set
\[
    \mathbb{P}(\mathcal{F})=\{\mathfrak{p}\in \mathbb{P}(\mathbb{N}):\mathfrak{p}|n\text{ for some } n\in \mathcal{F}\}.
\]
For $n\in \mathcal{F}\setminus\{1\}$, consider prime factorisations of the form
\[
    n=\prod_{\mathfrak{p}\in \mathbb{P}(\mathcal{F})} \mathfrak{p}^{\nu_{\mathfrak{p}}(n)}
\]
where the exponent $\nu_{\mathfrak{p}}(n)=0$ for $\mathfrak{p}$  not dividing $n$. Let
\[
    d=    \# \mathbb{P}(\mathcal{F})<\infty.
    \]
Order the elements of $\mathbb{P}(\mathcal{F})$ in an increasing manner so that
\[
    \mathbb{P}(\mathcal{F})=\{\mathfrak{p}_1<\ldots<\mathfrak{p}_d\}.
\]
Then
\[
    n=\mathfrak{p}_1^{\nu_{\mathfrak{p}_1}\!(n)}\cdots
    \mathfrak{p}_d^{\nu_{\mathfrak{p}_d}\!(n)}
    \qquad \forall n\in \mathcal{F}\setminus\{1\}.
\]
Below we allow $d=0$, for  $\mathcal{F}=\{1\}$ and $\mathbb{P}(\mathcal{F})=\varnothing$.

Let $\{c_n\}_{n\in \mathcal{F}}\subset \mathbb{C}$. 
The finite Dirichlet series 
\[
     m(z)=\sum_{n\in \mathcal{F}}\frac{c_n}{n^z} \qquad z\in \mathbb{C}
\]
is naturally identified with a polynomial in $d$ variables, as follows. Without ambiguity $p(w)=c_1$ whenever $d=0$. 
For $d\geq 1$ consider the $d$-dimensional polydisk,
\[
 \mathbb{D}^d=\{(w_1,  ...,w_d): \max_{j \in \mathbb{N}}\vert w_j \vert < 1\}
\]
with its distinguished boundary
\[
 \mathbb{T}^d=\{(w_1, ...,w_d): \vert w_j \vert=1 \  \forall j=1,\ldots,d\}.
\]
 Let
\[
     p(w)=p(w_1,\ldots,w_d)=\sum_{n\in \mathcal{F}} c_{n}  w_{1}^{\nu_{\mathfrak{p}_1}\!(n)}\cdots w_{d}^{\nu_{\mathfrak{p}_d}\!(n)}.
\]
Here and elsewhere,  $w=(w_1,\ldots,w_d)$.
Then
\[
     m(z)=\sum_{n\in \mathcal{F}}\frac{c_n}{\mathfrak{p}_1^{z\nu_{\mathfrak{p}_1}\!(n)}\cdots \mathfrak{p}_d^{z\nu_{\mathfrak{p}_d}\!(n)}}=p(\mathfrak{p}_1^{-z},\ldots, \mathfrak{p}_d^{-z}).
\]
Moreover, by the maximum principle,
\[
     \sup_{\Re(z)>0}|m(z)|=\sup_{w\in \mathbb{D}^d}|p(w)|=
     \max_{w\in \mathbb{T}^d}|p(w)|
\]
and also
\[
     \inf_{\Re(z)>0}|m(z)|=\inf_{w\in \mathbb{D}^d}|p(w)|=
     \min_{w\in \mathbb{T}^d}|p(w)|.
\]

\begin{mylem} \label{cri3}
Let $f:\mathbb{R}\longrightarrow \mathbb{C}$ be a 2-periodic function. Assume that the Fourier coefficients of $f$ are such that 
\[
|\hat{f}(j)|\leq \phi_j \qquad \forall j\in \mathbb{N}
\]
for a sequence $\{\phi_j\}_{j=1}^\infty\in \ell^1(\mathbb{N})$ and 
let $\varphi=\sum_{j=1}^\infty \phi_j$. Let $\mathcal{F}\subset \mathbb{N}$ be a finite set such that $1\in\mathcal{F}$ and let
\[
     \mu=\min_{w\in \mathbb{T}^d}\left|\sum_{j\in \mathcal{F}}\hat{f}(j) w_{1}^{\nu_{\mathfrak{p}_1}\!(j)}\cdots w_{d}^{\nu_{\mathfrak{p}_d}\!(j)}  \right|.
\]
Let $k\in \mathbb{N}$. If 
\begin{equation} \label{invereasy}
\sum_{j\in \mathcal{F}\setminus\{1\}} 
|\hat{f}({j})|<\hat{f}(1)\end{equation}
and
\begin{equation} \label{twomodes}
     \mu-\varphi+\sum_{j\in \mathcal{F}} |\hat{f}({j})| +\sum_{j=1}^k \left(\phi_j-|\hat{f}(j)|\right)>0,
\end{equation}
then $\mathfrak{E}_f\approxident \mathfrak{E}_{\operatorname{s}}$.
\end{mylem}
\begin{proof}
Decompose
\[
     m_f(z)=m(z)+v(z)
\qquad \text{where} \qquad
    m(z)=\sum_{j\in \mathcal{F}} \frac{\hat{f}({j})}{j^z}.    
\]
Let
\[
p(w)=\sum_{j\in \mathcal{F}} \hat{f}({j})  w_{1}^{\nu_{\mathfrak{p}_1}\!(j)}\cdots w_{d}^{\nu_{\mathfrak{p}_d}\!(j)}.
\]
From \eqref{invereasy} it follows that $p(w)$ has all its zeros in the complement of $\overline{\mathbb{D}^d}$ and so $\mu>0$. Now
\begin{align*}
   \sup_{\Re(z)>0}|v(z)|&\leq
   \sum_{j\in\mathbb{N}} |\hat{f}(j)|-\sum_{j\in\mathcal{F}} |\hat{f}(j)| \\ & < \sum_{j=1}^k |\hat{f}(j)|+\sum_{j=k+1}^\infty \phi_j -\sum_{j\in\mathcal{F}} |\hat{f}(j)| \\
   &= \sum_{j=1}^k (|\hat{f}(j)|-\phi_j)+\varphi -\sum_{j\in\mathcal{F}} |\hat{f}(j)|<\mu
\end{align*}
where the last inequality is implied by \eqref{twomodes}.
Hence
\[
     \sigma=\sup_{\Re(z)>0}\frac{|v(z)|}{|m(z)|}<1
\]
and so
\[
   \inf_{\Re(z)>0}|m_f(z)|= \inf_{\Re(z)>0}\left|m(z)\left(1+\frac{v(z)}{m(z)}\right)\right|\geq \mu(1-\sigma)>0.
\]
\end{proof}

Consider the following consequence of this lemma. Let $\mathfrak{p}\in\mathbb{P}(\mathbb{N})$. Assume that
\begin{equation} \label{crihome0}
    \hat{f}({\mathfrak{p}^2})>0 \quad \text{and} \quad \hat{f}({\mathfrak{p}^2})+|\hat{f}(\mathfrak{p})|<\hat{f}(1).
    \end{equation}
Either of the following two hypotheses ensure that $\mathfrak{E}_f\approxident \mathfrak{E}_{\operatorname{s}}$, \cite[corollaries~4.3 and 4.4]{BL2014}.
\begin{enumerate}
\item $|\hat{f}(\mathfrak{p})|\left[\hat{f}({\mathfrak{p}^2})+\hat{f}(1)\right]\geq 4\hat{f}({\mathfrak{p}^2})\hat{f}(1)$ and
\begin{equation}   \label{crihome1} 
\sum_{j \in \mathbb{N}\setminus \{1, \mathfrak{p}^2\}
  } \vert \hat{f}(j) \vert <\hat{f}(1)+\hat{f}({\mathfrak{p}^2}),
\end{equation}
\item $|\hat{f}(\mathfrak{p})|\left[\hat{f}({\mathfrak{p}^2})+\hat{f}(1)\right]<4\hat{f}({\mathfrak{p}^2})\hat{f}(1)$ and
\begin{equation}   \label{crihome2} 
\sum_{j \in \mathbb{N}\setminus \{1, \mathfrak{p},\mathfrak{p}^2\}
  } \vert \hat{f}(j) \vert <\left[\hat{f}(1)-\hat{f}({\mathfrak{p}^2})\right]\sqrt{1-\frac{[\hat{f}(\mathfrak{p})]^2}{4\hat{f}(1)\hat{f}({\mathfrak{p}^2})}}.
\end{equation}
\end{enumerate}
For a proof of this, set $d=1$, $k=0$ and $\mathcal{F}=\{1, \mathfrak{p},
\mathfrak{p}^2\}$. According to \cite[Lemma~4.1]{BL2014}, if \eqref{crihome0} holds true, then $\mu>0$. Moreover, 
    \begin{equation} \label{forlaterlemma}
     |\hat{f}(\mathfrak{p})|\left[\hat{f}({\mathfrak{p}^2})+\hat{f}(1)\right]\geq 4\hat{f}({\mathfrak{p}^2})\hat{f}(1)\ \Rightarrow \ \mu=\hat{f}(1)+\hat{f}({\mathfrak{p}^2})-|\hat{f}(\mathfrak{p})|
   \end{equation}
   and
    \[
     |\hat{f}(\mathfrak{p})|\left[\hat{f}({\mathfrak{p}^2})+\hat{f}(1)\right]< 4\hat{f}({\mathfrak{p}^2})\hat{f}(1)\ \Rightarrow \ \mu=
     \left[\hat{f}(1)-\hat{f}({\mathfrak{p}^2})\right]\sqrt{1-\frac{[\hat{f}({\mathfrak{p}})]^2}{4\hat{f}(1)\hat{f}({\mathfrak{p}^2})}}
   \]  
Therefore the hypotheses of Lemma~\ref{cri3} are satisfied whenever \eqref{crihome1} or \eqref{crihome2} hold. 

\begin{remark}
Most likely a version of Lemma~\ref{cri3} can be established for the Banach space setting $L^r(0,1)$, by following the ideas announced in the recent work \cite{Mit2017}. However, various details need to be carefully confirmed. 
\end{remark}


\section{Piecewise linear profiles} \label{section3}
Let $0<\alpha\leq \frac12$. Set
\begin{equation*} 
 g_\alpha(x)= 
\begin{cases}
 \frac{x}{\alpha} & 0 \leq x<\alpha\\
1 & \alpha \leq x<1-\alpha\\
\frac{1-x}{\alpha} & 1-\alpha \leq x\leq 1.
\end{cases}
\end{equation*}
Extend $g_\alpha$ to an odd function on $[-1,1]$ then to a $2$-periodic function on $\mathbb{R}$. It is known \cite[\S 5]{BBCDG2006}  that $\mathfrak{E}_{g_{\frac12}}\approxident \mathfrak{E}_{\operatorname{s}}$. This section addresses the full equivalence of $\mathfrak{E}_{g_{\alpha}}$ with the Fourier basis for $\alpha$ near $0$.  

\subsection{Fourier coefficients}
Since
\begin{align*}
 \widehat{g_\alpha}(j)&=\frac{2}{\alpha}\int_0^\alpha x \sin(j\pi x)\mathrm{d}x+\int_\alpha^{1-\alpha} \sin(j \pi x)\mathrm{d}x+\frac{2}{\alpha}\int_{1-\alpha}^1 (1-x)\sin(j \pi x)\mathrm{d}x,
\end{align*}
then
\begin{equation}\label{Scoeff}
{\widehat{g_\alpha}(j)=\begin{cases} 0 & \forall j\equiv_2 0 \\ 
\frac{4}{\alpha j^2 \pi^2}\sin(j \pi \alpha) & \forall j\equiv_2 1.\end{cases}}
\end{equation}

The proof of the next lemma follows a similar path as the argument described in 
\cite[p.49]{BushellEdmunds2012}. We include details.

\begin{mylem}   \label{lemcoe} For all $0<\alpha<\frac12$,
   \begin{align*}
       \sum_{j=1}^\infty \widehat{g_\alpha}({j})&=\frac{2}{\alpha}\int_0^\alpha \frac{x}{\sin(\pi x)}\mathrm{d}x+\frac{2}{\pi}\log \frac{1+\cos(\alpha \pi)}{\sin(\alpha \pi)}.   
       \end{align*}
\end{mylem}
\begin{proof}
 Let $r\in[0,1]$ and $\lambda(r)=\frac{2r}{1+r^2}$. Let
 \[
     \phi(r,\alpha)=2\int_{0}^{1/2}  \frac{g_\alpha(x)\lambda(r) \sin (\pi x)}{1-\lambda^2(r) \cos^2(\pi x)}\ \mathrm{d} x .
 \]
 Then 
\[
    \phi(1,\alpha)=2\int_{0}^{1/2} \frac{g_\alpha(x)}{\sin(\pi x)} \ \mathrm{d}x,
\] 
where the integral is finite because $g_\alpha(x)$ is linear near $x=0$. Now
 \[
      2\sum_{k=0}^\infty
     r^{2k+1} \sin((2k+1)\pi x)=\frac{\lambda(r) \sin (\pi x)}{1-\lambda^2(r) \cos^2(\pi x)}
     \qquad \forall r\in[0,1)
 \]
 where the series on the left hand side is absolutely convergent.
Then, by the dominated convergence theorem,
 \[
     \phi(r,\alpha)=4\int_0^{1/2} g_{\alpha}(x) \sum_{k=0}^\infty
     r^{2k+1} \sin((2k+1)\pi x) \ \mathrm{d}x=
     \sum_{k=0}^\infty r^{2k+1} \widehat{g_\alpha}({2k+1}) 
      \]
 for all $r\in[0,1)$.

Now, from \eqref{Scoeff} it follows that
$
    \sum_{k=0}^\infty |\widehat{g_\alpha}({2k+1})|<\infty.
$ 
 Then the series $\sum_{k=0}^\infty \widehat{g_\alpha}({2k+1})$ is  absolutely convergent. By virtue of Abel's limit theorem, we have
\[
        \sum_{k=0}^\infty \widehat{g_\alpha}({2k+1})=
        \lim_{r\to 1}\phi(r,\alpha)=\phi(1,\alpha)=
        2\int_{0}^{1/2} \frac{g_\alpha(x)}{\sin(\pi x)} \ \mathrm{d}x.
\]
Hence
\begin{align*}
 \sum_{k=0}^\infty \widehat{g_\alpha}({2k+1})
&=\frac{2}{\alpha}\int_0^\alpha \frac{x}{\sin(\pi x)}\mathrm{d}x+2\int_\alpha^\frac12 \frac{\mathrm{d}x}{\sin(\pi x)}\\
&=\frac{2}{\alpha}\int_0^\alpha \frac{x}{\sin(\pi x)}\mathrm{d}x+\frac{2}{\pi}\log \frac{1+\cos(\alpha \pi)}{\sin(\alpha \pi)}.
\end{align*}
\end{proof}

\subsection{Basis properties of $\mathfrak{E}_{g_\alpha}$}
Since 
\[
    \sum_{j=3}^\infty |\widehat{g_\alpha}(j)|<\frac{4}{\alpha \pi^2}\sum_{j=1}^\infty \frac{1}{(2j+1)^2}=
    \frac{4}{\alpha \pi^2}\left(\frac{\pi^2}{8}-1    \right) 
\]
and
\[
    \widehat{g_\alpha}(1)=\frac{4\sin(\pi \alpha)}{\alpha \pi^2},
\]
then
\[
    \sum_{j=3}^\infty |\widehat{g_\alpha}(j)|<\widehat{g_\alpha}(1)
\]
whenever
\[
      \sin(\pi \alpha)>\frac{\pi^2}{8}-1.
\]
As $\sin(\pi \alpha)$ is increasing in $\alpha\in (0,\frac{1}{2})$, then $\mathfrak{E}_{g_\alpha}\approxident \mathfrak{E}_{\operatorname{s}}$ for
all $\alpha\in (\alpha_0,\frac12]$ where \[\alpha_0=\frac{1}{\pi}\arcsin\left(\frac{\pi^2}{8}-1\right)\approx 0.0750835.\] 
In the following lemma, this threshold is moved towards $\alpha=0$ by quite a significant margin.

\begin{mylem} \label{leforak}
If  $0<\alpha<\frac{1}{2}$ is such that
\begin{equation}  \label{eqleforak} \tag{$\ast$}
     2\sin(\pi\alpha)+\sum_{j=0}^k \frac{1-|\sin((2j+1)\pi \alpha)|}{(2j+1)^2}
     >
     \frac{\pi^2}{8}
\end{equation}
for some $k\in \mathbb{N}$,  then $\mathfrak{E}_{g_\alpha}\approxident \mathfrak{E}_{\operatorname{s}}$. 
\end{mylem}
\begin{proof}
This is a consequence of Lemma~\ref{cri3} with
\[
     \phi_j=\begin{cases}0& j\equiv_2 0 \\ \frac{4}{\alpha j^2 \pi^2}& j\equiv_2 1 \end{cases}, \qquad m(z)=\widehat{g_\alpha}(1)=\frac{4\sin(\alpha \pi)}{\alpha\pi^2}
\]
and $\mathcal{F}=\{1\}$. In this case
\[
    \varphi=\frac{1}{2\alpha} \quad \text{and}\quad   \mu=\frac{4\sin(\alpha \pi)}{\alpha\pi^2}.
\]
Notice that the hypothesis \eqref{invereasy} is trivial and that \eqref{eqleforak} is a re-arrangement of \eqref{twomodes}.
\end{proof}

\begin{figure}
\centerline{\includegraphics[width=.5\textwidth]{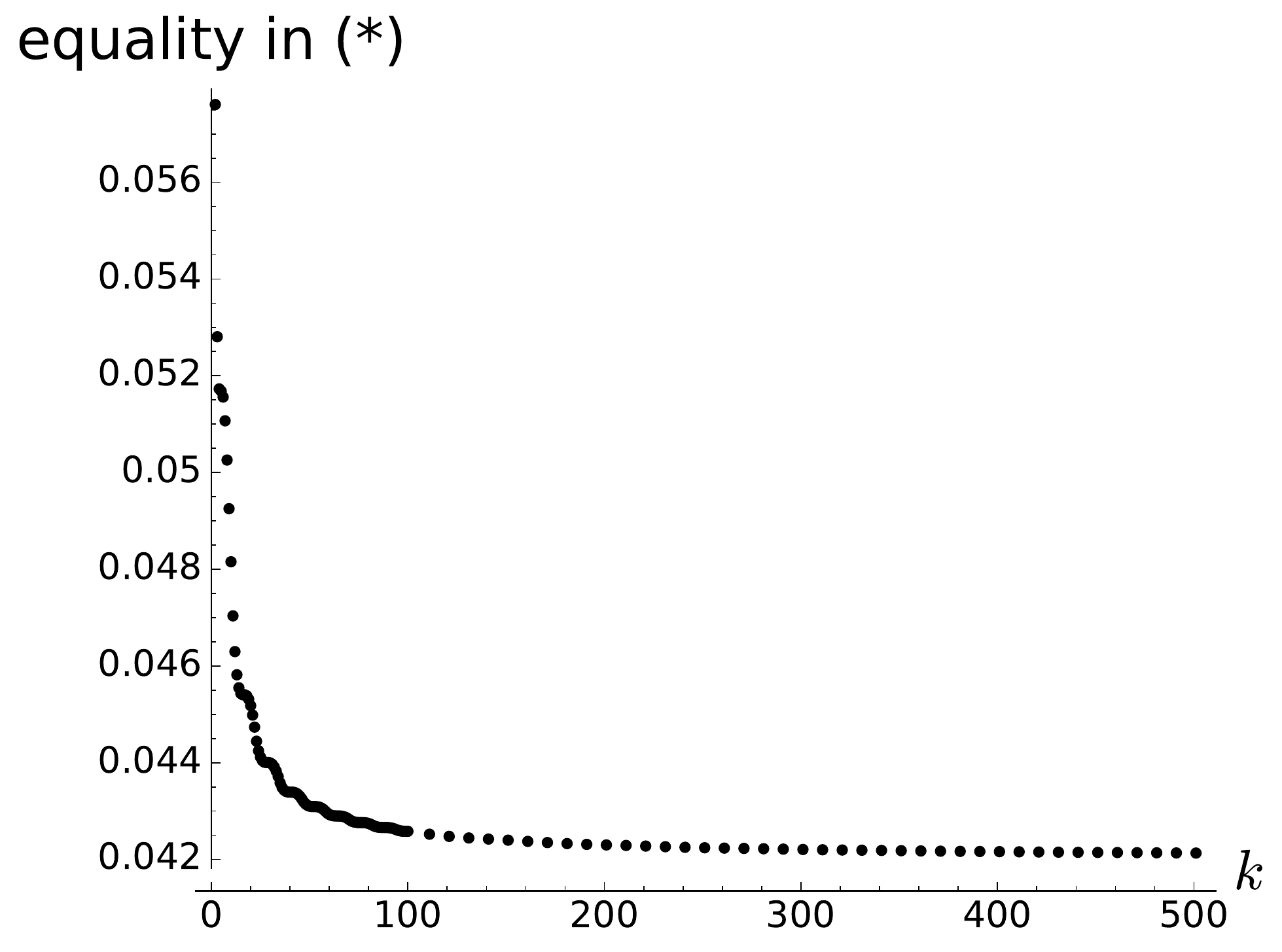}
\includegraphics[width=.5\textwidth]{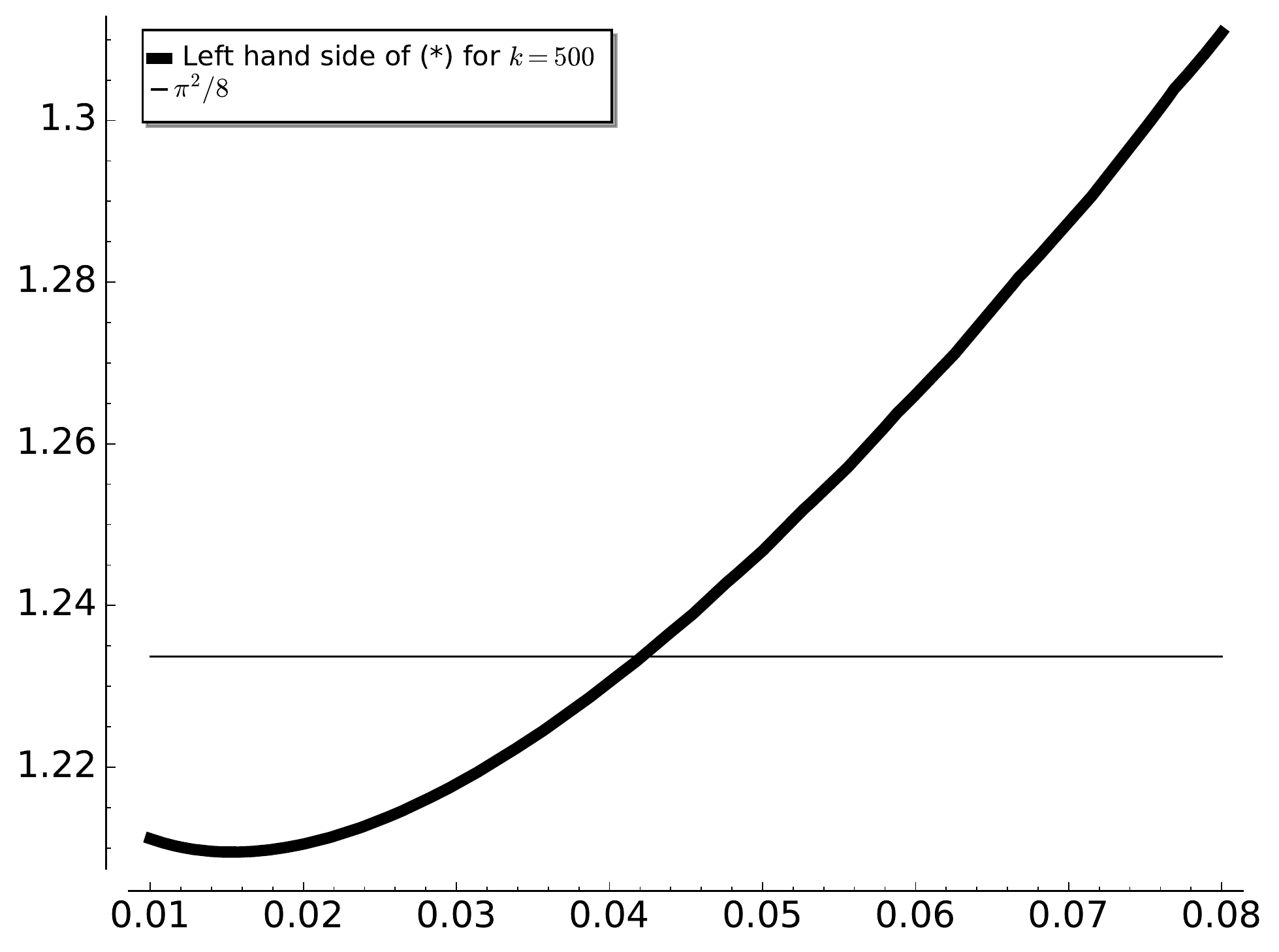}}
\caption{Left: values of $\alpha$ (vertical axis) where equality is attained in \eqref{eqleforak} of Lemma~\ref{leforak} for increasing $k$. The step size is 1 for $k\leq 100$ and it is 10 for $100<k\leq 500$. Right: left hand side of \eqref{eqleforak} for $k=500$ for $\alpha\in[\frac{1}{100},\frac{8}{100}]$.  \label{figlemmaforak}}
\end{figure}

As $k$ increases, equality in \eqref{eqleforak} is achieved for smaller values of $\alpha$. For small values of $k$, the behaviour of the root in terms of $k$ is oscillatory and quite complicated but it eventually stabilises as $k$ increases. See the left of Figure~\ref{figlemmaforak}. For $k=500$ a numerical approximation of the solution of the equation 
\[
2\sin(\pi\alpha)+\sum_{j=0}^{500} \frac{1-|\sin((2j+1)\pi \alpha)|}{(2j+1)^2}
     =
     \frac{\pi^2}{8}
\]
is \[\alpha_1\approx 0.0421317.\] 
The right side of Figure~\ref{figlemmaforak} shows graphically that \eqref{eqleforak} holds true for all
$\alpha\in(\alpha_1,\frac{8}{100}]$. An analytic cofirmation of this would be rather tedious and probably not worth pursuing.

Taking $k$ larger in Lemma~\ref{leforak} would not allow confirmation of full equivalence for $\alpha$ much closer to $0$. For $k$ beyond $500$, the tail of the summation would only contribute by a factor smaller than $10^{-3}$ on the left hand side of \eqref{eqleforak}. In turn 
\[
    \widehat{g_{\alpha_2}}(1)= \sum_{j=3}^\infty |\widehat{g_{\alpha_2}}(j)|
\]
where $\alpha_2< \alpha_1$ matches the first 2 significant figures. A numerical approximation of both sides of the following expression indicates that
\[
     \sum_{j=3}^{111}|\widehat{g_{\alpha_3}}(j)|>\widehat{g_{\alpha_3}}(1) \qquad \text{for} \qquad \alpha_3=0.04.
\]

Moreover, recall Lemma~\ref{lemcoe}. From the latter it follows that the identity
\[
 \sum_{j=0}^\infty \widehat{g_\alpha}({2j+1})=2\widehat{g_\alpha}(1)
\]
is satisfied for 
\[
 \pi^2\int_{0}^\alpha \frac{x}{\sin (\pi x)}\mathrm{d}x
 +\pi \alpha \log\frac{1+\cos(\alpha \pi)}{\sin(\alpha \pi)}=
 4\sin(\alpha \pi).
\]
A numerical solution to this is $\alpha=\alpha_4\approx 0.0318993$.  Eventually, for $\alpha\to 0$, 
\[
    \widehat{g_{\alpha}}(1)\leq \sum_{j=3}^\infty |\widehat{g_\alpha}(j)|.
\]
In order to detect further the confirmed threshold for basis of $\mathfrak{E}_{g_\alpha}$, it is necessary to take $d=2$, $\mathfrak{p}_1=3$ and $\mathfrak{p}_2=5$ as follows.

\begin{figure}
\centerline{\includegraphics[width=.6\textwidth]{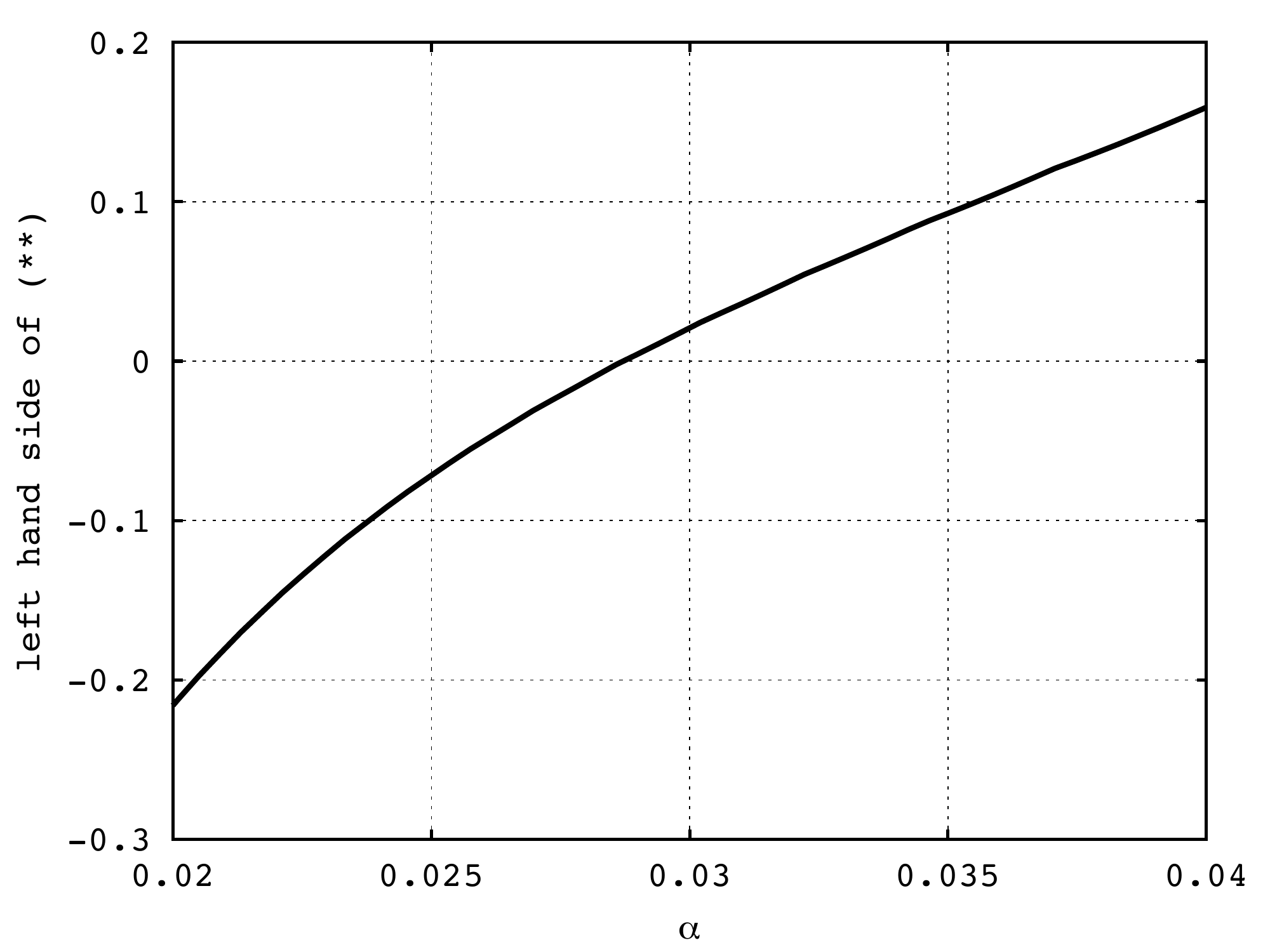}}
\caption{Left hand side of \eqref{eqleforak2} for $k=50$ and $\alpha\in[\frac{2}{100},\frac{4}{100}]$ \label{figlemmaforak2}.}
\end{figure}

\begin{mylem} \label{leforak2}
Let
\[
    \begin{aligned}\hat{\mu}(\alpha)=\min_{(x,y)\in[-\pi,\pi]^2} & \big[\widehat{g_\alpha}(1)+\widehat{g_\alpha}(3)\cos(x)+\widehat{g_\alpha}(9)\cos(2x)+ \\
& \qquad \qquad \qquad \widehat{g_\alpha}(5)\cos(y)+  \widehat{g_\alpha}({25})\cos(2y)\big]^2+ \\
    &\big[\widehat{g_\alpha}(3)\sin(x)+\widehat{g_\alpha}(9)\sin(2x)+ \\
&  \qquad \qquad \qquad \widehat{g_\alpha}(5)\sin(y)+\widehat{g_\alpha}({25})\sin(2y)\big]^2.\end{aligned}
\]
If  $0<\alpha<\frac{1}{25}$ is such that
\begin{equation}  \label{eqleforak2} \tag{$\ast\ast$}
\sqrt{\hat{\mu}(\alpha)}-
     \frac{\pi^2}{8}+
\sum_{j\in\{1,3,5,9,25\}}\frac{\sin(j\pi \alpha)}{j^2}+\sum_{j=0}^k \frac{1-|\sin((2j+1)\pi \alpha)|}{(2j+1)^2}
     >0
\end{equation}
for some $k\in \mathbb{N}$,  then $\mathfrak{E}_{g_\alpha} \approxident \mathfrak{E}_{\operatorname{s}}$. 
\end{mylem}
\begin{proof}
In Lemma~\ref{cri3} once again set
\[
     \phi_j=\begin{cases}0& j\equiv_2 0 \\ \frac{4}{\alpha j^2 \pi^2}& j\equiv_2 1 \end{cases}
     \] 
but then set $d=2$ and 
\[    
     p(w)=p(w_1,w_2)=\widehat{g_\alpha}(1)+\widehat{g_\alpha}(3)w_1+\widehat{g_\alpha}(5)w_2+\widehat{g_\alpha}(9)w_1^2+\widehat{g_\alpha}({25})w_2^2.
\]
Note that \eqref{invereasy} holds true for all $0<\alpha<\frac{1}{25}$. Indeed $\widehat{g_\alpha}(1)>0$. Also,
\begin{align*}
     \frac{\cos(3\pi \alpha)}{3}&+\frac{\cos(5\pi \alpha)}{5}+\frac{\cos(9\pi \alpha)}{9}+\frac{\cos(25\pi \alpha)}{25} \\ &<
    \left(\frac{1}{3}+\frac{1}{5}+\frac{1}{9}+\frac{1}{25}\right)\cos(\pi \alpha) \\
&< \cos(\pi \alpha)
\end{align*}
for all $0<\alpha<\frac{1}{25}$. 
Then,
\[
     \frac{\sin(3\pi \alpha)}{9}+\frac{\sin(5\pi \alpha)}{25}+\frac{\sin(9\pi \alpha)}{81}+\frac{\sin(25\pi \alpha)}{625}<\sin(\pi \alpha)
\] 
for all such  $\alpha$, ensuring the validity of \eqref{invereasy}.

Now, parametrise $(w_1,w_2)\in \mathbb{T}^2$ by means of
\[ 
w_1=\cos(x)+i\sin(x) \qquad \text{and} \qquad
w_2=\cos(y)+i\sin(y)
\]
where $(x,y)\in[-\pi,\pi]^2$. Then  
\[
    \hat{\mu}(\alpha)=\min_{w\in \mathbb{T}^2}|p(w)|^2.
\]
Re-arranging the condition \eqref{twomodes} leads to the condition
\eqref{eqleforak2}.
\end{proof}

For $k=50$ equality in \eqref{eqleforak2} is achieved for $\alpha=\alpha_5\approx 0.0287740$. The picture in Figure~\ref{figlemmaforak2} shows that \eqref{eqleforak2} holds true for all $\alpha\in(\alpha_5,\frac{1}{25}]$.


\section{Continuously differentiable profiles} \label{section4}
Let $0<\beta<\frac12$. Set
\[
   h_\beta(x)=\begin{cases}
      \left(\frac{x}{\beta}+1 \right)^2\left(1-\frac{x}{2\beta}\right)-1 & 0\leq x < \beta \\
      1 & {\beta} \leq x\leq \frac{1}{2}
   \end{cases}
\]
Extend $h_\beta$ to $[0,1]$ by reflection at $\frac12$, then to an odd function in $[-1,1]$ and then to a 2-periodic function on $\mathbb{R}$. The derivative $h_\beta'(x)$ is continuous on $\mathbb{R}$. Moreover, $h_\beta(x)\to \operatorname{J}(x)$ as $\beta\to 0$. This section examines the full equivalence of the family $\mathfrak{E}_{h_\beta}$ and the Fourier basis for $\beta$ near $0$.

\subsection{The Fourier coefficients}

The Fourier coefficients of $h_\beta$ are
\begin{equation}\label{Scoeff2}
{\widehat{h_\beta}(j)=\begin{cases} 0 & \forall j\equiv_2 0 
\\ 
 \frac{12}{j^3\pi^3\beta^2}\left[\frac{\sin(j\pi \beta)}{j\pi \beta}-  \cos (j\pi \beta)\right] & \forall j\equiv_2 1.\end{cases}}
\end{equation}
Put
\[
   \phi_j=\begin{cases}  0 & j\equiv_2 0 \\
   \frac{12}{\pi^3\beta^2}\left[\frac{1}{j^4\pi\beta}+\frac{1}{j^3}\right] &  j\equiv_2 1.  \end{cases}
\]
Then 
\[
     |\widehat{h_\beta}(j)|\leq \phi_j.   
\]
Also note that
\[
    \widehat{h_\beta}(j)=\mathrm{O}(j^{-3}) \qquad \forall \beta>0.
\]

\subsection{Basis properties} \label{sec3}

We firstly consider the simplest case.

\begin{mylem}   \label{exa3basic}
If $0<\beta<\frac{1}{2}$ is such that
\begin{equation} \label{exa3first}
        \left[\frac{\pi^3}{96}-\frac{1}{\pi}\right]\beta^{-1}+\frac{7}{8}\zeta(3)-1<\frac{\sin(\pi \beta)}{\pi\beta}-\cos (\pi \beta),
\end{equation}
then $\mathfrak{E}_{h_\beta}\approxident \mathfrak{E}_{\operatorname{s}}$.
\end{mylem}
\begin{proof}
Let
\[
\gamma_1=\sum_{\substack{j=3 \\ j\not\equiv_2 0}}\frac{1}{j^3}=\frac{7}{8}\zeta(3)-1
\qquad
\text{and} \qquad 
    \gamma_2=\sum_{\substack{j=3 \\ j\not\equiv_2 0}}\frac{1}{j^4}=
    \frac{\pi^4}{96} - 1   .
\]
Then
\begin{align*}
   \sum_{j=3}^\infty |\widehat{h_\beta}(j)|& \leq \sum_{j=3}^\infty\frac{12}{j^3\pi^3 \beta^2}
   \left[\frac{1}{j\pi\beta}+1\right] \\
   &= \frac{12}{\pi^3 \beta^2}\left[\frac{\gamma_2}{\pi\beta}+\gamma_1   \right]
\end{align*}
and\[
     \widehat{h_\beta}(1)=\frac{12}{\pi^3 \beta^2}\left[\frac{\sin(\pi \beta)}{\pi\beta} -\cos (\pi \beta)    \right]>0.
\]

Hence the condition \eqref{exa3first} implies
\[
     \sum_{j=3}^\infty |\widehat{h_\beta}(j)|<\widehat{h_\beta}(1).
\]
\end{proof}

In this lemma, the left hand side of \eqref{exa3first} is decreasing in $\beta$ and it has a singularity $+\infty$ as $\beta\to 0$. On the other hand, the right hand side is increasing from the value $0$ at $\beta=0$. Moreover, \eqref{exa3first} holds true for $\beta=\frac12$ and equality is achieved for 
\[
\beta=\beta_0\approx 0.159059.
\]
Thus, full equivalence is ensured for all $\beta\in(\beta_0,\frac12]$.

\begin{remark}
The Fourier coefficients of $h_\beta$ also satisfy the inequality
\[
       |\widehat{h_\beta}(j)|\leq \frac{24}{j^3\pi^3 \beta^2 } \qquad \forall j\in \mathbb{N}.
\]
Then the condition
\[
     \frac{7}{4}\zeta(3)-2<\frac{\sin(\pi \beta)}{\pi\beta}-\cos (\pi \beta)
\]
also yields $\mathfrak{E}_{h_\beta}\approxident \mathfrak{E}_{\operatorname{s}}$. As it turns, this other condition only holds true for $\beta\in(\tilde{\beta}_0,\frac12]$ where $\tilde{\beta}_0\approx 0.180340$. 
\end{remark}

 \begin{figure}
\centerline{\includegraphics[width=.8\textwidth]{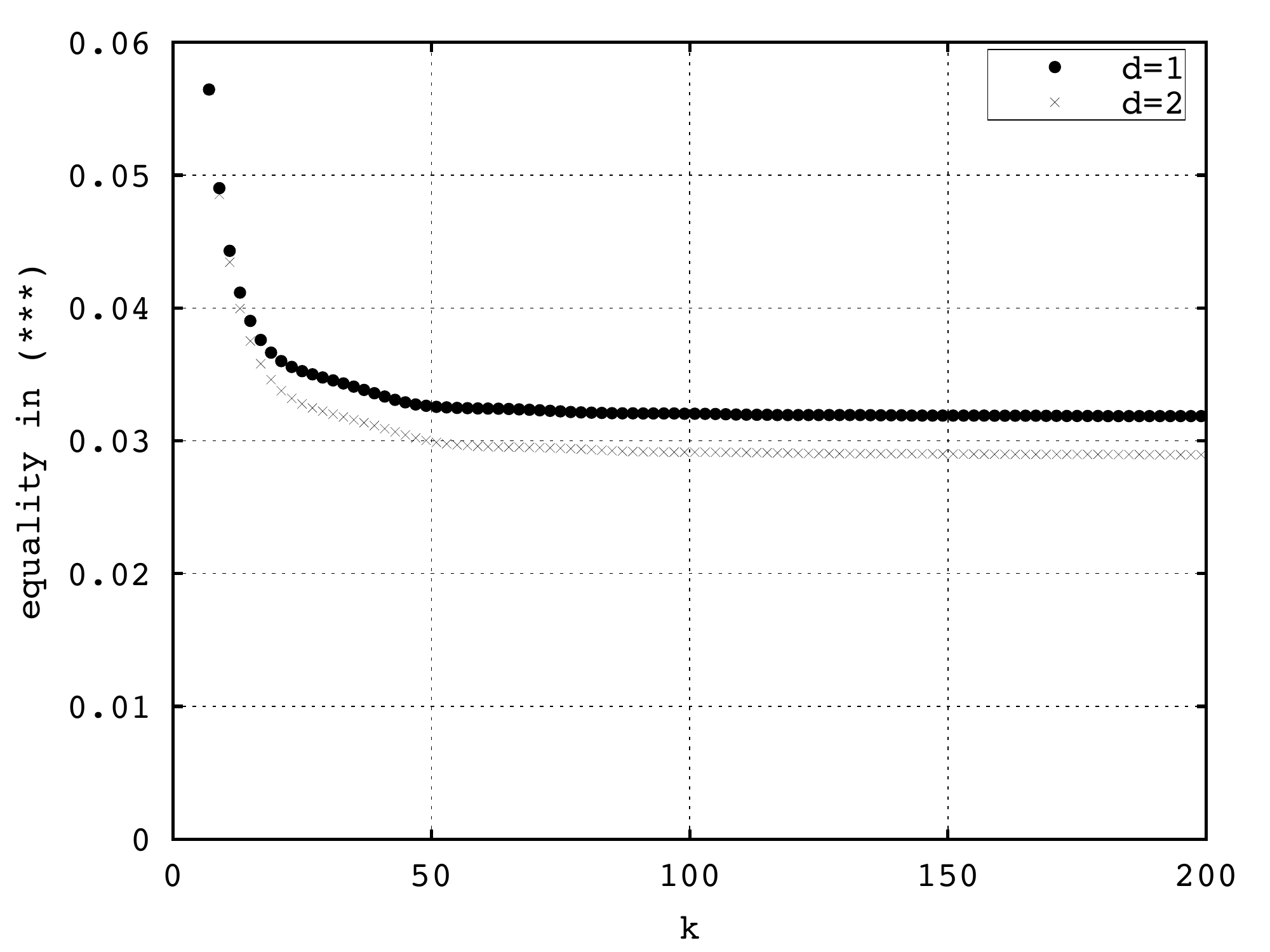}}
\centerline{\includegraphics[width=.5\textwidth]{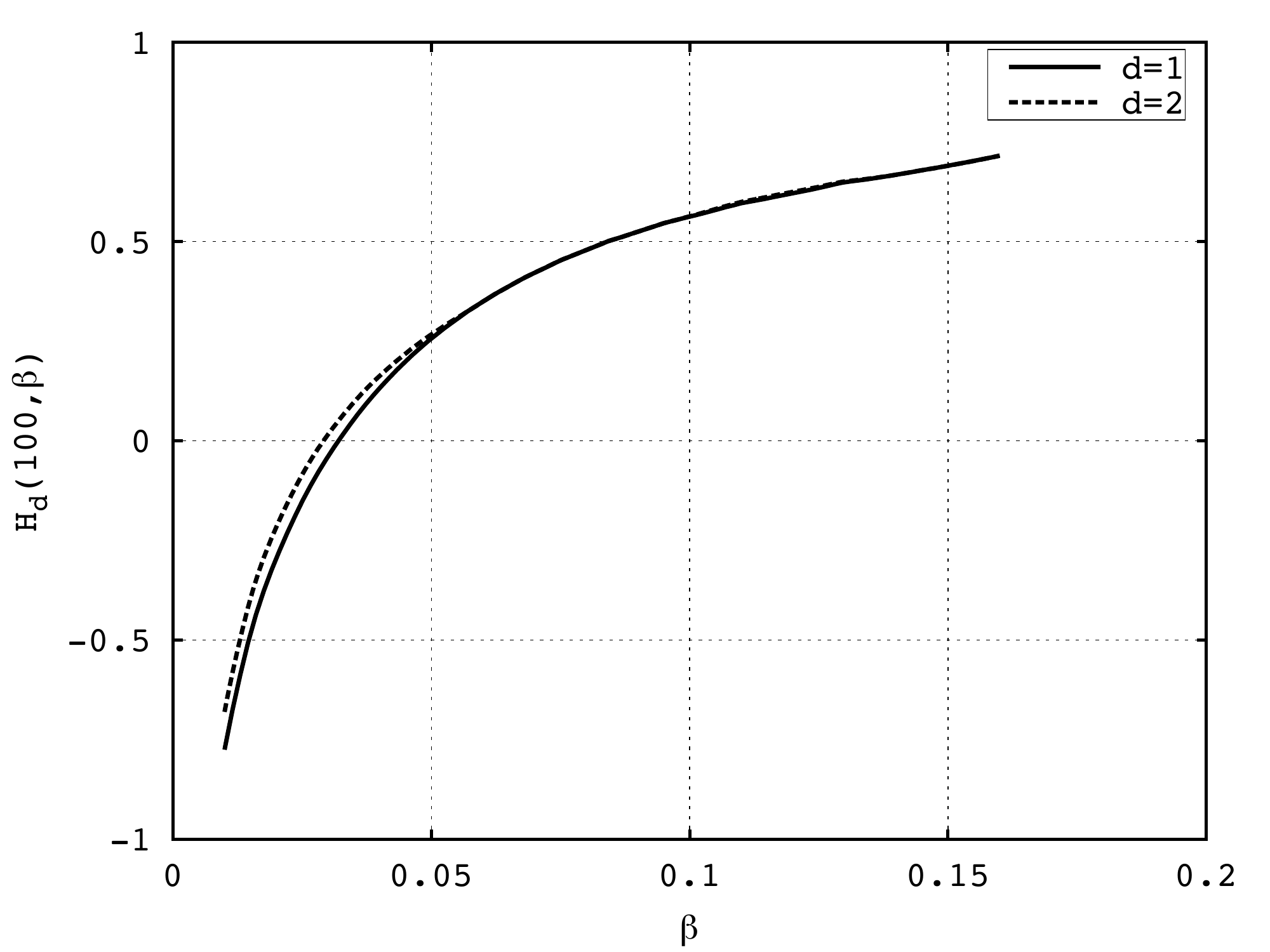}
\includegraphics[width=.5\textwidth]{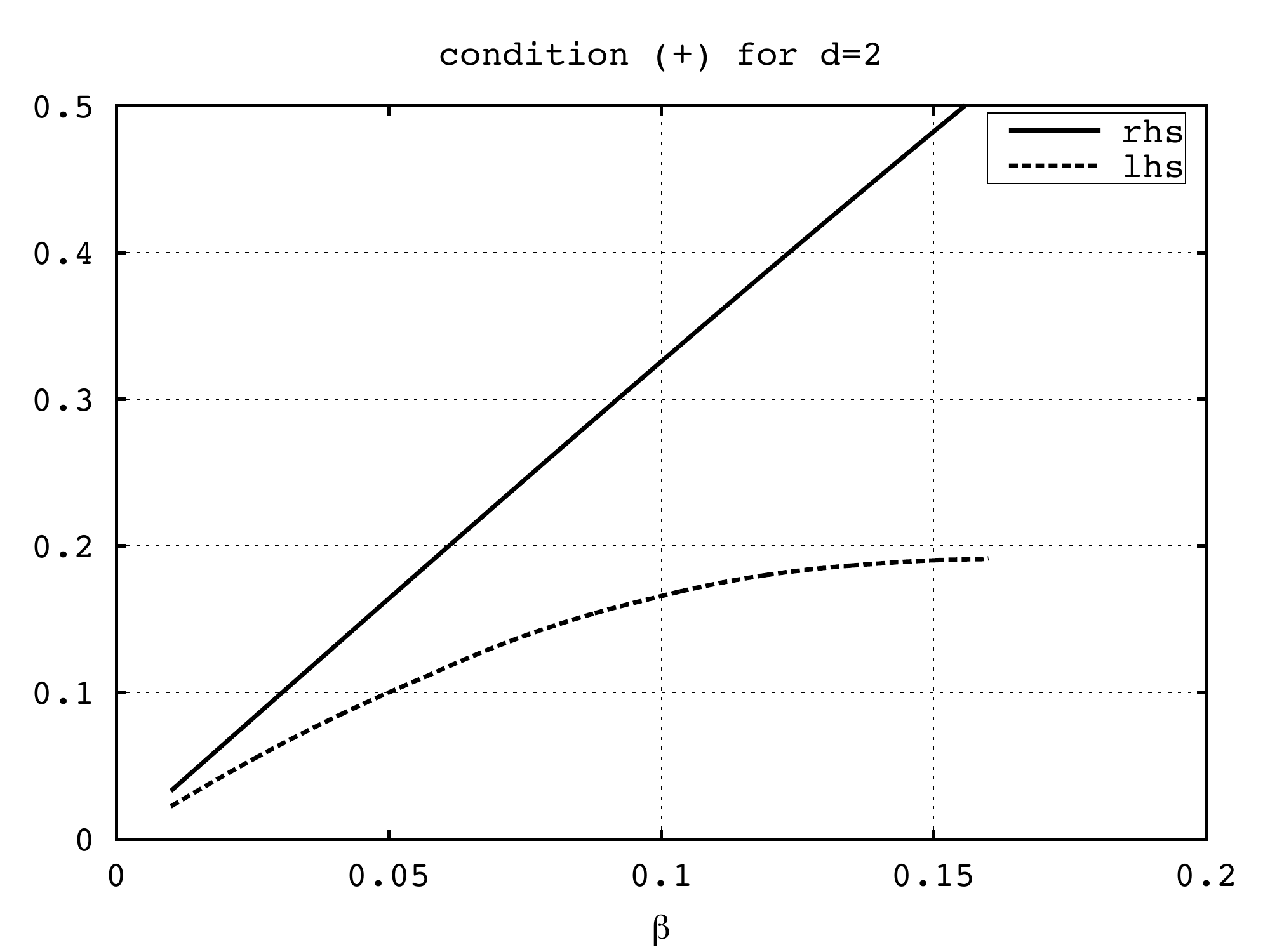}
}
\caption{top: value of $\beta$ (vertical axis) where $H_d(k,\beta)=0$ for $k$ increasing. Bottom left: plot of $H_d(100,\beta)$ for $\beta\in[\frac{1}{100},\frac{16}{100}]$.
Bottom right: condition \eqref{invereasybj} in the case $d=2$ for $\beta\in[\frac{2}{100},\frac{16}{100}]$. The picture indicates that this condition is clearly satisfied in the shown range.   \label{figforbk}}
\end{figure}

For different parameters in Lemma~\ref{cri3}, the threshold $\beta_0$ is improved by a significant margin. Put
\[
     \varphi=\tilde{\varphi}(\beta)=\frac{12}{\pi^3\beta^2}
     \left[\frac{1}{\pi \beta} \sum_{\substack{j=1 \\ j\equiv_2 1}}^\infty \frac{1}{j^4} + \sum_{\substack{j=1 \\ j\equiv_2 1}}^\infty \frac{1}{j^3} \right]=\frac{12}{\pi^3\beta^2}\left[\frac{\pi^3}{96\beta}+\frac{7}{8}\zeta(3)   \right].
\]
Fix $\mathcal{F}=\{1,\mathfrak{p}_1,\mathfrak{p}_1^2,\ldots,\mathfrak{p}_d,\mathfrak{p}_d^2\}$.
Set
\[
     \begin{aligned}
     H_d(k,\beta)&= \tilde{\mu}(\beta)-\tilde{\varphi}(\beta)+
     \sum_{j\in \mathcal{F}}|\widehat{h_\beta}(j)|+
     \sum_{j=1}^k\left(\phi_j- |\widehat{h_\beta}(j)| \right)
     \end{aligned}
\]
where
\[
     \tilde{\mu}(\beta)=\min\left\{\left| \widehat{h_\beta}(1)+\sum_{j=1}^d 
      \widehat{h_\beta}({\mathfrak{p}_j})w_j+\widehat{h_\beta}({\mathfrak{p}_j^2})w_j^2 \right|:
      w\in \mathbb{T}^d    \right\}.
\]
That is, $H_d(k,\beta)$ is the left hand side of \eqref{twomodes} for the periodic function $h_\beta$. The following lemma is a consequence of Lemma~\ref{cri3}.

\begin{mylem} \label{ex2multimodes}
Let $0<\beta<\frac12$ and $k,\,d\in \mathbb{N}$. If
\begin{equation} \label{invereasybj} \tag{$+$}
\sum_{j\in\{\mathfrak{p}_1,\mathfrak{p}_1^2,\ldots,\mathfrak{p}_d,\mathfrak{p}_d^2\}} \left|\frac{\sin(j\pi \beta)}{j^4\pi \beta^2}-\frac{\cos(j\pi \beta)}{j^3\beta}\right|<
\frac{\sin(\pi \beta)}{\pi\beta^2}-\frac{\cos (\pi \beta)}{\beta}
\end{equation}
and
\begin{equation} \label{exa2modes} \tag{$\ast\ast\ast$}
     H_d(k,\beta)>0 ,
\end{equation}
then 
$\mathfrak{E}_{h_\beta}\approxident \mathfrak{E}_{\operatorname{s}}$.
\end{mylem}

Figure~\ref{figforbk} shows graphical confirmations of the hypotheses  of Lemma~\ref{ex2multimodes} for the two cases, $d=1$ and $d=2$. Approximation of solutions of the equations 
\[
    H_1(100,\beta)=0 \qquad \text{and} \qquad H_2(100,\beta)=0
\]
are
\[
     \beta_1\approx 0.0320481 \qquad \text{and} \qquad
     \beta_2\approx 0.0291447,
\]
respectively. The graphs on the bottom indicate that
both hypotheses of Lemma~\ref{ex2multimodes} for $d=2$ hold in the range $\beta\in(\beta_2,\beta_0]$. Therefore in that range also $\mathfrak{E}_{h_\beta}\approxident \mathfrak{E}_{\operatorname{s}}$ is expected.


\section{$p$-sine function profiles} \label{section5}

The basis properties of the next benchmark example have been examined in the series of papers \cite{BBCDG2006, BushellEdmunds2012, EdmundsGurkaLang2012, BL2014}. 

Let $p>1$ and  $p'=\frac{p}{p-1}$. Denote by $B(a,b)$ and $\mathcal{I}(a,b;t)$, the beta and the incomplete beta functions in their usual parameters \cite[8.391 $\text{and}$ 8.392]{Grad2007}.  Let $F_p:[0,1]\longrightarrow [0,\pi_p/2]$ be given by
\[
F_p(y)=\int_0^y\frac{\mathrm{d}t}{(1-t^p)^{\frac{1}{p}}}  = \frac{\pi_p}{2}\mathcal{I}\left(\frac1p, \frac{1}{p'}; y^p\right)
\]
where
\[
    \pi_p=2F_p(1)=\frac{2B\Big(\frac{1}{p},\frac{1}{p'}\Big)}{p}=\frac{2\pi}{p\sin(\frac{\pi}{p})}.
\]
The $p$-sine function $\sin_p:\mathbb{R}\longrightarrow [-1,1]$ is defined as the inverse function
\[
     \sin_p(x)=F^{-1}_p(x) \qquad x\in\left[0,\frac{\pi_p}{2}\right]
\]
extended by the rules
\[
    \sin_p(-x)=-\sin_p(x) \qquad \text{and} \qquad 
    \sin_p\left(\frac{\pi_p}{2}-x\right)=\sin_p\left(\frac{\pi_p}{2}+x\right),
\]
which make this function $2\pi_p$-periodic, differentiable, odd with respect to 0 and even with respect to $\pi_p/2$. Note that $\sin_2(x)=\sin(x)$ and $\pi_2=\pi$.

Let
\[
     \operatorname{s}_p(x)=\sin_p(\pi_p x). 
\]
In \cite{BBCDG2006} and \cite{BushellEdmunds2012} it was determined that $\mathfrak{E}_{\operatorname{s}_p}\approxident \mathfrak{E}_{\operatorname{s}}$ for all $p>p_1$
where $p_1\approx 1.19824$. This threshold was subsequently improved in \cite{BL2014} to $p>p_2\approx 1.04392$. As seen next, a suitable application of Lemma~\ref{cri3} lowers the range of full equivalence to the Fourier basis to a point closer to $p=1$ by a significant margin.

\subsection{Fourier coefficients}
Integration by parts and changing to $ t=\sin_p(\pi_p x)$ the variable of integration  yield
\begin{align*}
 \widehat{\mathrm{s}_p}(j)
&=4 \int_0^\frac{1}{2} \sin_p(\pi_p x) \sin(j \pi x) \mathrm{d}x\\
&=\frac{4}{j \pi}\int_0^\frac{1}{2} [\sin_p(\pi_p x)]' \cos(j \pi x) \mathrm{d}x\\
&=\frac{4}{j\pi}\int_0^1 \cos\left[\frac{j \pi}{\pi_p}F_p(t)\right]\mathrm{d}t.
\end{align*}
Hence  
\begin{equation}\label{ajnewform}
 {\widehat{\mathrm{s}_p}(j)=\begin{cases} 0 & \forall j\equiv_2 0 \\
\frac{4}{j\pi}\int_0^1 \cos\left[\frac{j \pi}{\pi_p}F_p(t)\right]\mathrm{d}t
& \forall j\equiv_2 1 
\end{cases}}
\end{equation}

The next inequality \cite[\S4 (4.3)]{BushellEdmunds2012} will be employed in several places below,
\begin{align}\label{aupperbound}
|\widehat{\mathrm{s}_p}(j)| <\frac{4 \pi_p}{j^2\pi^2} \qquad \forall j \geq 1,\,p>1.
\end{align}
It is known \cite{BushellEdmunds2012} that
\[
     \widehat{\mathrm{s}_{p_3}}(1)=\sum_{j=3}^\infty \widehat{\mathrm{s}_{p_3}}(j)
\]
for $p_3\approx  1.04399$. Note that $p_2<p_3$.
 Various new technical points about $\widehat{\mathrm{s}_p}(j)$  are included in the Appendix.


\subsection{Basis properties}
 
In \cite{BL2014} the threshold $p_2$ mentioned above for full equivalence to the Fourier basis was obtained as a consequence of a statement \cite[Proposition~7.1]{BL2014} very similar to the following lemma.

\begin{mylem} \label{forsine}
Let $k\geq 9$ and $1<p<\frac{12}{11}$. Suppose that
\begin{enumerate}
\item $|\widehat{\mathrm{s}_p}(3)|+\widehat{\mathrm{s}_p}(9)<\widehat{\mathrm{s}_p}(1)$,
\item $|\widehat{\mathrm{s}_p}(3)|[\widehat{\mathrm{s}_p}(1)+\widehat{\mathrm{s}_p}(9)]  \geq 4 \widehat{\mathrm{s}_p}(9) \widehat{\mathrm{s}_p}(1)$.
\end{enumerate}
If
\begin{align}\label{consr2}
\frac{\pi_p}{2}-\frac{4\pi_p}{\pi^2}\sum_{\substack{
   j =1\\
   j\equiv_2 1
  }}^k \frac{1}{j^2}
 < \widehat{\mathrm{s}_p}(1)+\widehat{\mathrm{s}_p}(9) -\sum_{\substack{j=3 \\ j\not=9}}^k  |\widehat{\mathrm{s}_p}(j)| ,
\end{align}
then $\mathfrak{E}_{\operatorname{s}_p}\approxident \mathfrak{E}_{\operatorname{s}}$.
\end{mylem}
\begin{proof}
This lemma is a consequence of Lemma~\ref{cri3}. Put $\mathcal{F}=\{1,3,9\}$,
\begin{equation} \label{phijforsinep}
  \phi_j=\begin{cases}  0 & j\equiv_2 0 \\
     \frac{4\pi_p}{j^2\pi^2}&  j\equiv_2 1  \end{cases}, \qquad \qquad \varphi=\frac{\pi_p}{2}
\end{equation}
and recall \eqref{aupperbound}. According to \cite[Lemma~6.1]{BL2014}, $\widehat{\mathrm{s}_p}(9)>0$ for all $1<p<\frac{12}{11}$. Thus, the condition~1 implies \eqref{invereasy}. Moreover, both condition~1 and 2, imply
\[
     \mu=\widehat{\mathrm{s}_p}(1)+\widehat{\mathrm{s}_p}(9)-|\widehat{\mathrm{s}_p}(3)|,
\] 
see \eqref{forlaterlemma}.
With this data, \eqref{twomodes} reduces exactly to \eqref{consr2}.
\end{proof}

\begin{remark}
There are two minor differences between the Proposition~7.1 of \cite{BL2014} and Lemma~\ref{forsine} above. 
In the former it was additionally required that all the Fourier coefficients $\widehat{\mathrm{s}_p}(j)\geq 0$ for $1\leq j\leq k$. On the other hand, 
Lemma~\ref{forsine} includes the extra condition~2, which ensures the hypothesis \eqref{invereasy} automatically. See also the assumption \eqref{crihome0}. This is sufficient, but not necessary, for $\mu>0$. 
\end{remark}

According to the calculations performed in \cite[p21]{BL2014}, the condition~2 of Lemma~\ref{forsine} is satisfied for any $p \in (p_4, \frac{12}{11})$ where $p_4 \approx 1.03854$. Numerical verification indicates that condition~1 holds true for $p\in(1.01,1.1)$. See also \cite[Lemma~4.2]{BL2014}.  Estimating  both sides of \eqref{consr2} 
for different values of $k$, indicates that equality occurs in this identity for $p\approx 1.03876>p_4$ when $k=61$ and $p\approx 1.03852<p_4$ for $k=63$. From this information, it follows that an application of Lemma~\ref{forsine} only extends the threshold for full equivalence up to $p_4$ and not beyond that point.

As seen next, an analogue to Lemma~\ref{ex2multimodes} in this context moves the threshold further towards $p=1$.  Put $\mathcal{F}=\{1,\mathfrak{p}_1,\mathfrak{p}_1^2,\ldots,\mathfrak{p}_d,\mathfrak{p}_d^2\}$. Consider the same choice \eqref{phijforsinep}. Set
\[
     \begin{aligned}
     J_d(k,p)&= \bar{\mu}(p)-\frac{\pi_p}{2}+
     \sum_{j\in \mathcal{F}}|\widehat{\mathrm{s}_p}(j)|+
     \sum_{j=1}^k\left(\phi_j-|\widehat{\mathrm{s}_p}(j)| \right)
     \end{aligned}
\]
where
\[
     \bar{\mu}(p)=\min\left\{\left| \widehat{\mathrm{s}_p}(1)+\sum_{j=1}^d 
      \widehat{\mathrm{s}_p}({\mathfrak{p}_j})w_j+\widehat{\mathrm{s}_p}({\mathfrak{p}_j^2})w_j^2 \right|:
      w\in \mathbb{T}^d    \right\}.
\]
Recall \eqref{twomodes}.

\begin{mylem} \label{ex3multimodes}
Let $p>1$ and $k,\,d\in \mathbb{N}$. If
\begin{equation} \label{invereasycj} \tag{$++$}
\sum_{j\in\{\mathfrak{p}_1,\mathfrak{p}_1^2,\ldots,\mathfrak{p}_d,\mathfrak{p}_d^2\}} |\widehat{\mathrm{s}_p}(j)|<
\widehat{\mathrm{s}_p}(1)
\end{equation}
and
\begin{equation} \label{exa3modes} \tag{${\ast\ast}{\ast\ast}$}
     J_d(k,p)>0 ,
\end{equation}
then 
$\mathfrak{E}_{\operatorname{s}_p}\approxident \mathfrak{E}_{\operatorname{s}}$.
\end{mylem}

\begin{figure}
\centerline{\includegraphics[width=.5\textwidth]{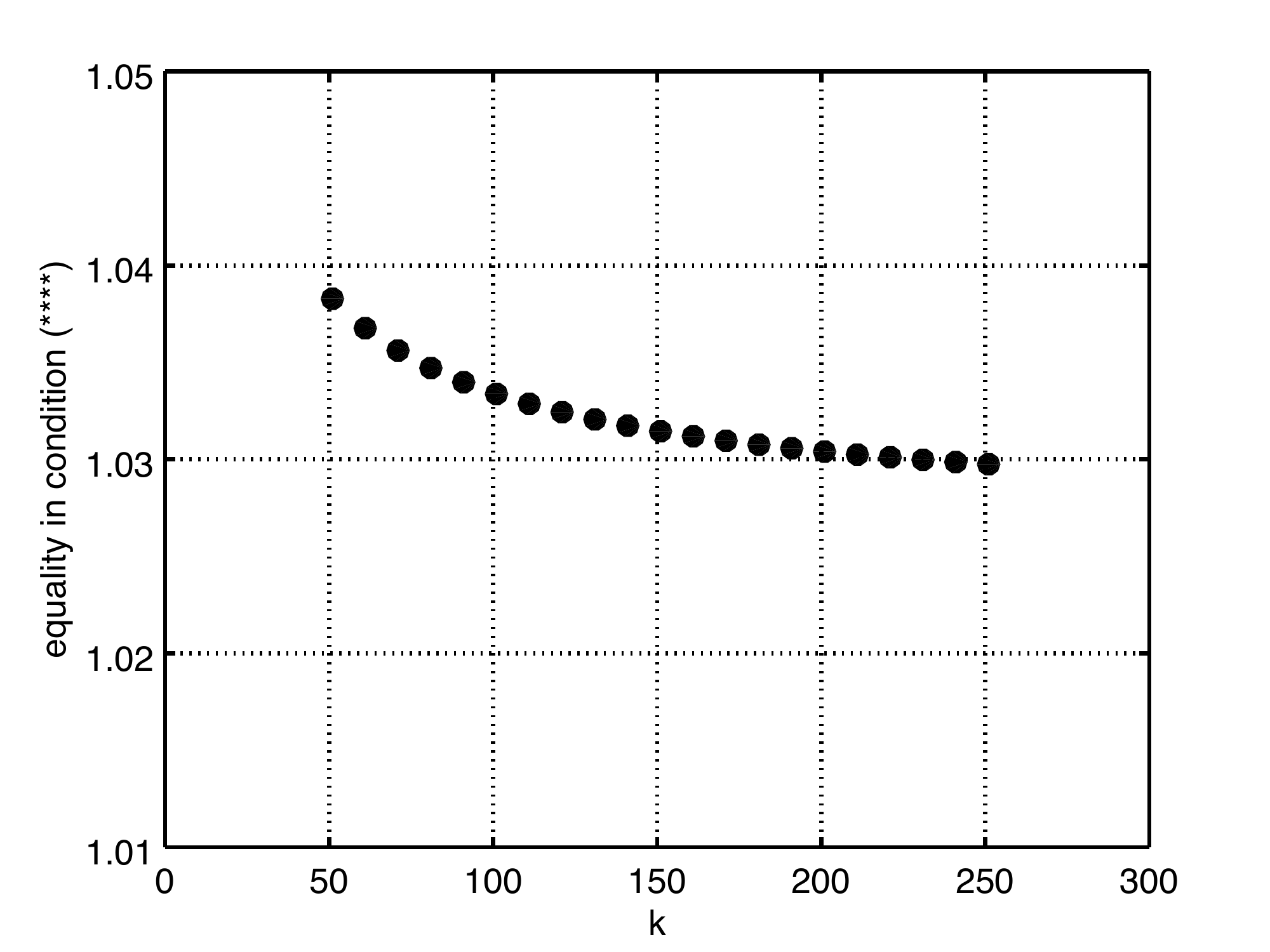}
\includegraphics[width=.5\textwidth]{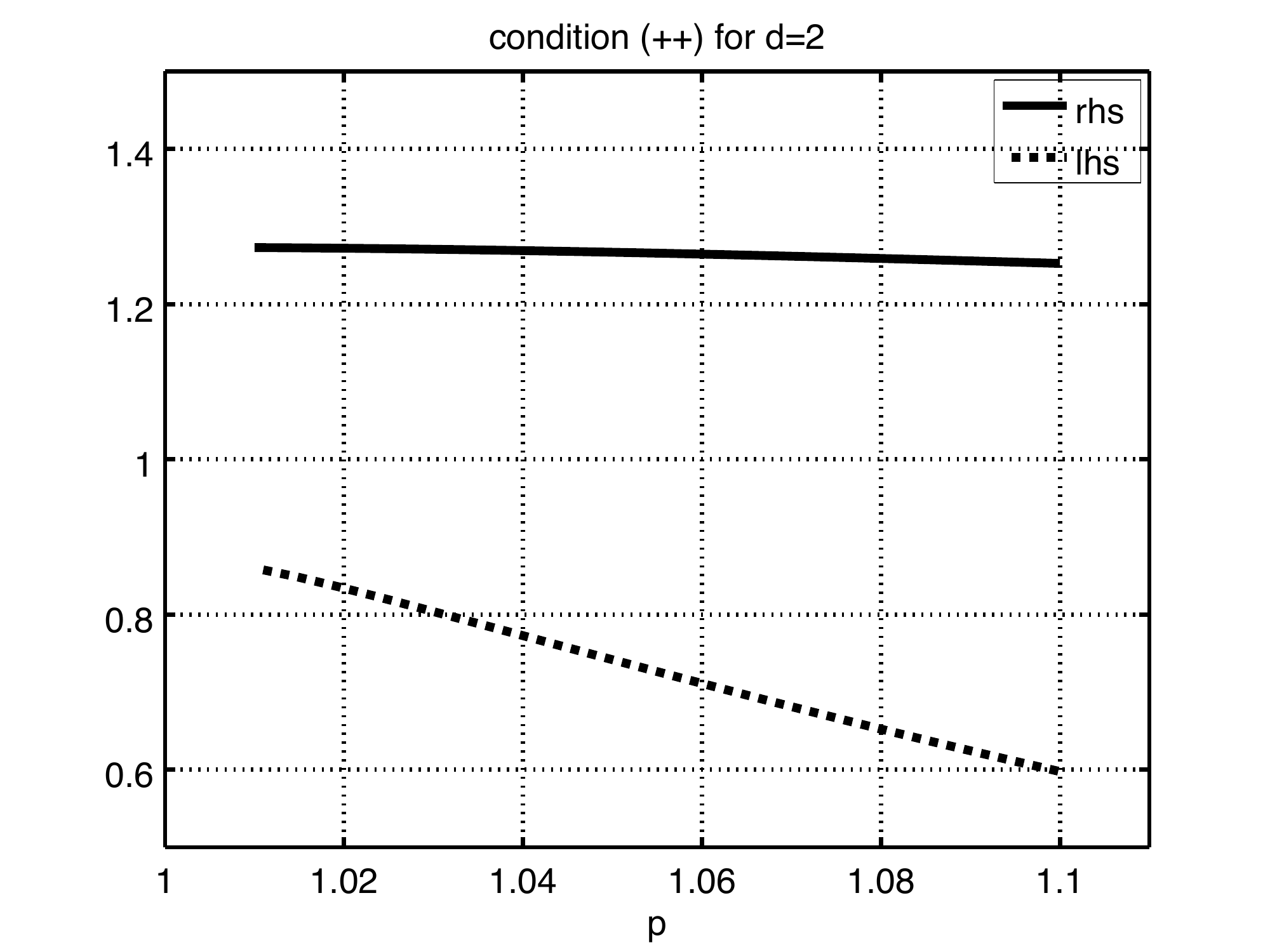}}
\caption{Graphs corresponding to Lemma~\ref{ex3multimodes}. Left: value of $p$ (vertical axis) where $J_2(k,p)=0$ for $k$ increasing.   
Right: both sides of the condition \eqref{invereasycj} for $d=2$.
\label{figforck}}
\end{figure}

Take $d=2$ in this lemma, so $\mathfrak{p}_1=3$ and $\mathfrak{p}_2=5$. As $k$ increases, equality in \eqref{exa3modes} is achieved for decreasing values of $p$. See Figure~\ref{figforck} (left). For $k=251$, a numerical approximation of the equation yields equality for $p=p_5\approx 1.02975$ and $J_2(251,p)$ is an increasing function for $p\in[1.01,1.1]$.  
Moreover, the condition \eqref{invereasycj} holds true for $d=2$ in this range. See Figure~\ref{figforck} (right). This indicates that
 $\mathfrak{E}_{\operatorname{s}_p}\approxident \mathfrak{E}_{\operatorname{s}}$
for all $p>p_5$.  


 
\appendix

\section{Towards analytical bounds for $\widehat{\mathrm{s}_p}(j)$}

The estimation of the value $p_5$ obtained  above is the best threshold for full equivalence to the Fourier basis of the $p$-sine functions that we can report at present time. In this appendix we include various  estimates for the Fourier coefficients of $\mathrm{s}_p$ which might be used for analytic confirmation of this threshold.

\subsection{Properties of the inverse $\sin_p$ function} \label{theappen}

We begin by recalling the following fundamental property established in \cite[Corollary~4.4]{BushellEdmunds2012}. 
Let
\[
    \mathcal{I}_p(x)=\frac{2}{\pi_p}F_p(x)=\mathcal{I}\left(\frac1p, \frac{1}{p'}; y^p\right).
\]
For any $p, q \in (1, \infty)$ such that $p < q$,
\begin{align}\label{incbetaincp}
     1<\frac{\mathcal{I}_q(y)}{\mathcal{I}_p(y)}<\frac{\pi_p}{\pi_q}.
\end{align}

\begin{mylem}\label{incbetaconvex}
For $p>1$ fixed. The function $\mathcal{I}_p(y)$ is monotonically increasing and convex in $y \in [0, 1]$.
\end{mylem}
\begin{proof}
Indeed
\[
\frac{\mathrm{d}}{\mathrm{d}y}\mathcal{I}_p(y)=\frac{2}{\pi_p}(1-y^p)^{-\frac{1}{p}}>0
\]
and
\[
\frac{\mathrm{d}^2}{\mathrm{d}y^2}\mathcal{I}_p(y)=\frac{2}{\pi_p}y^{p-1}(1-y^p)^{-\frac1p -1}>0.
\]
\end{proof}
\begin{mylem}\label{segment}
Let $k\in \mathbb{N}$ be fixed. Let $0 \leq y< x \leq 1$ be such that $\cos\big(\frac{k \pi}{2}\mathcal{I}_p(u) \big)$ is decreasing for all $u \in [y, x]$.
Then,
\begin{align*}
 \int_{y}^{x} \cos\left(\frac{k\pi}{2}\mathcal{I}_p(u)\right) \mathrm{d}u>\mathcal{I}_{k, p}(y, x)
 \end{align*}
 where
 \begin{align*}
 \mathcal{I}_{k, p}(y, x):=\frac{2}{k\pi}(x-y)\frac{\sin\big(\frac{k\pi}{2}\mathcal{I}_p(x)\big)-\sin\big(\frac{k\pi}{2}\mathcal{I}_p(y)\big)}{\mathcal{I}_p({x})-\mathcal{I}_p(y)}.
 \end{align*}
\end{mylem}
\begin{proof}
The chord of $\mathcal{I}_p(u)$ with endpoints $y$ and ${x}$ is given by
\begin{align}\label{chord}
 f(u)=\frac{\mathcal{I}_p({x})-\mathcal{I}_p(y)}{{x}-y}(u-{x})+\mathcal{I}_p({x}).
\end{align}
Lemma~\ref{incbetaconvex} implies, $\mathcal{I}_p(u)<f(u)$ for any $p>1$ and $u \in (0, 1)$. Hence and by virtue of the hypothesis,
\begin{align*}
 \int_{y}^{{x}} \cos\left(\frac{k\pi}{2}\mathcal{I}_p(u)\right) \mathrm{d}u >  \int_{y}^{{x}} \cos\left(\frac{k\pi}{2}f(u)\right) \mathrm{d}u=\mathcal{I}_{k, p}(y, {x}).
\end{align*}
 \end{proof}
For $0<s<{t}<1$ consider the function
\begin{equation} \label{functionG}
G(s, {t}):=\frac{(1-s^p)^\frac{1}{p}\mathcal{I}_p(s)-(1-{t}^p)^\frac{1}{p}\mathcal{I}_p({t})+\frac{2}{\pi_p}({t}-s)}{(1-s^p)^\frac{1}{p}-(1-{t}^p)^\frac{1}{p}}.
\end{equation}

\begin{mylem}\label{tangentab}
Let $k\in \mathbb{N}$ be fixed. Let $0 \leq s< {t} \leq 1$ be such that $\cos\left(\frac{k \pi}{2}\mathcal{I}_p(u) \right)$ is increasing for all $u \in [s, {t}]$. Then
\begin{align*}
 \int_{s}^{{t}} \cos\left(\frac{k\pi}{2}\mathcal{I}_p(u)\right) \mathrm{d}u &> \mathcal{J}_{k, p}^{(1)}(s,{t})+\mathcal{J}_{k, p}^{(2)}(s, {t})
\end{align*}
where
\begin{align*}
\mathcal{J}_{k, p}^{(1)}(s, {t})&:=\frac{\pi_p}{k \pi}(1-s^p)^\frac{1}{p}\left[\sin\left(\frac{k\pi}{2}G(s, {t})\right)-\sin\left(\frac{k\pi}{2}\mathcal{I}_p(s)\right)\right],\\
\mathcal{J}_{k, p}^{(2)}(s, {t})&:=\frac{\pi_p}{k \pi}(1-{t}^p)^\frac{1}{p}\left[\sin\left(\frac{k\pi}{2}\mathcal{I}_p({t})\right)-\sin\left(\frac{k\pi}{2}G(s, {t})\right)\right].
\end{align*}
\end{mylem}
\begin{proof}
The tangent to the curve $\mathcal{I}_p(u)$ at any $u=s$ is given by
\begin{equation}\label{tangent}
\gamma_s(u)=\frac{2}{\pi_p (1-s^p)^\frac{1}{p}}(u-s)+\mathcal{I}_p(s).
\end{equation}
By virtue of Lemma~\ref{incbetaconvex}, for any $p>1$ and $u\in [0,1]$, we have $\mathcal{I}_p(u)>\gamma_s(u)$. The intersection point $y$ of the tangents to $\mathcal{I}_p(u)$ at $s$ and ${t}$ is then given by
\begin{align*}
 y=\frac{\frac{\pi_p}{2}(1-{t}^p)^\frac{1}{p}(1-s^p)^\frac{1}{p}[\mathcal{I}_p(s)-\mathcal{I}_p({t})]+t(1-s^p)^\frac{1}{p}-s(1-{t}^p)^\frac{1}{p}}{(1-s^p)^\frac{1}{p}-(1-{t}^p)^\frac{1}{p}}.
\end{align*}
Moreover, $\gamma_s(y)=\gamma_{{t}}(y)=G(s, {t})$. \\
Then, because of the hypothesis,
\begin{align*}
 \int_s^{{t}} \cos\left(\frac{k\pi}{2}\mathcal{I}_p(u)\right) \mathrm{d}x &> \int_{s}^{y} \cos\left(\frac{k\pi}{2}\gamma_s(u)\right)\mathrm{d}u+\int_{y}^{{t}}\cos\left(\frac{k\pi}{2}\gamma_{{t}}(u)\right) \mathrm{d}u,
\end{align*}
where
\begin{align*}
\int_{s}^{y} \cos\left(\frac{k\pi}{2}\gamma_s(u)\right)\mathrm{d}u=\mathcal{J}_{k, p}^{(1)}(s, {t})
\end{align*}
and
\begin{align*}
\int_{y}^{{t}}\cos\left(\frac{k\pi}{2}\gamma_{{t}}(u)\right) \mathrm{d}u=\mathcal{J}_{k, p}^{(2)}(s, {t}).
\end{align*}
\end{proof}


\subsection{Towards analytic estimates for $\widehat{\mathrm{s}_p}(k)$ when $k\equiv_4 3$} \label{forckong3}
Let $k=4j-1$ for $j \in \mathbb{N}$. The integrand $\cos\left[\frac{k\pi}{2}\mathcal{I}_p(u)\right]$ in \eqref{ajnewform} for $u\in[0, 1]$ is monotonically decreasing in $j$ disjoint segments \[[\tilde{y}_i,\tilde{x}_i]\qquad \qquad i=1,\ldots,j\] and it is monotonically increasing in $j$ disjoint segments \[
[\tilde{s}_i,\tilde{t}_i] \qquad \qquad i=1,\ldots,j,\] 
so that
\[
      [0,1]=\left(\bigcup_{i=1}^{j} [\tilde{y}_i,\tilde{x}_i]  \right)
\cup
\left(\bigcup_{i=1}^{j} [\tilde{s}_i,\tilde{t}_i]  \right)
\]
where $\tilde{y}_1=0$, $\tilde{t}_{j}=1$, $\tilde{s}_i=\tilde{x}_i$ and $\tilde{y}_{i+1}=\tilde{t}_{i}$. The minimum turning points are  
such that
\[
      \mathcal{I}_p(\tilde{x}_i)=\frac{4m-2}{k}\qquad \text{for }m=1,\ldots,j
\]
and the maximum turning points are such that
\[
    \mathcal{I}_p(\tilde{t}_i)=\frac{4m}{k}\qquad \text{for }m=1,\ldots,j-1.
\]

We partition each one of these segments into sets of quadrature points as follows. Let $\{m_{i}^-\}_{i=1}^{j}\subset \mathbb{N}$ and $\{m_{i}^+\}_{i=1}^{j}\subset \mathbb{N}\setminus\{1\}$. Set 
\begin{align*}
&x_0=\tilde{y}_1=0,\quad x_{m_1^-}=\tilde{x}_1, \quad  x_{1+m_{1}^-}=\tilde{y}_2,\\ 
&x_{\sum_{\ell=1}^i m_\ell^-} =\tilde{x}_i,\quad x_{1+\sum_{\ell=1}^i m_\ell^-}=\tilde{y}_{i+1},  & (i=2,\ldots,j)\\
&t_1=\tilde{s}_1,\quad t_{m_1^+}=\tilde{t}_1,\quad t_{1+m_{1}^+}=\tilde{s}_2,\\
& t_{\sum_{\ell=1}^i m_\ell^+} =\tilde{t}_i,\quad  t_{1+\sum_{\ell=1}^i m_\ell^+}=\tilde{s}_{i+1}, & (i=2,\ldots,j)\\ 
& t_{\sum_{\ell=1}^j m_\ell^+} =\tilde{t}_j=1.
\end{align*}
We consider increasing sequences 
\begin{align*}
&0\leq \cdots<x_{m-1}<x_{m}<\cdots< 1 &(m=1,\ldots,\sum_{\ell=1}^{j} m_\ell^-) \\
&0< \cdots<t_{m-1}<t_{m}<\cdots \leq1 &(m=2,\ldots,\sum_{\ell=1}^{j} m_\ell^+)
\end{align*}
such that
\begin{gather*}
   \left\{x_{1+\sum_{\ell=1}^{i-1} m_\ell^-}<\cdots<
   x_{\sum_{\ell=1}^{i} m_\ell^-}   \right\} \subset
[\tilde{y}_{i},\tilde{x}_i] \qquad \text{and} \\
\left\{t_{1+\sum_{\ell=1}^{i-1} m_\ell^+}<\cdots<
   t_{\sum_{\ell=1}^{i} m_\ell^+}    \right\} \subset
[\tilde{s}_{i},\tilde{t}_i].
\end{gather*}

\begin{mylem}\label{a13lowerb}
Let $p>1$ and $k=4j-1$ where $j\in \mathbb{N}$. For $k>3$
\begin{align*}
 \widehat{\mathrm{s}_p}(k)>\frac{4}{k\pi} &\left[\sum_{m=1}^{m_1^-} \mathcal{I}_{k,p}(x_{m-1}, x_m)+\sum_{\ell=1}^{j-1} \sum_{m=\sum_{i=1}^{\ell}m_i^-+2}^{\sum_{i=1}^{\ell+1} m_{i}^-} \mathcal{I}_{k, p}(x_{m-1}, x_m) \right. \\  &\left.+\sum_{m=2}^{m_1^+} \Big(\mathcal{J}_{k, p}^{(1)}(t_{m-1}, t_m)+ \mathcal{J}_{k, p}^{(2)}(t_{m-1}, t_m)\Big)\right. \\ & \left. +\sum_{\ell=1}^{j-2}\sum_{m=\sum_{i=1}^{\ell}m_i^++2}^{\sum_{i=1}^{\ell+1}m_{i}^+}
\Big(\mathcal{J}_{k, p}^{(1)}(t_{m-1}, t_m)+\mathcal{J}_{k, p}^{(2)}(t_{m-1}, t_m)\Big)\right.\\ & \left. +\sum_{m=\sum_{i=1}^{j-1}m_i^++2}^{\sum_{i=1}^{j}m_{i}^+-1}
\Big(\mathcal{J}_{k, p}^{(1)}(t_{m-1}, t_m)+\mathcal{J}_{k, p}^{(2)}(t_{m-1}, t_m)\Big)\right.\\ &\left. +\mathcal{J}_{k, p}^{(1)}(t_{\sum_{\ell=1}^{j}m_{\ell}^+-1}, 1) \right].
\end{align*}
\end{mylem}
\begin{proof}
The proof follows from the properties established in lemmas~\ref{segment} and \ref{tangentab} by taking the endpoints as follows. 
\begin{itemize}
\item $y=x_{m-1}$ and $x=x_{m}$ in the former case. For $m=1, \ldots, m_1^-$ and $m=\sum_{i=1}^{\ell}m_i^-+~2, \dots, \sum_{i=1}^{\ell+1} m_{i}^-$ when $\ell=1, \ldots, {j-1}$. 
\item $s=t_{m-1}$ and $t=t_{m}$ in the latter case. For $m=2, \ldots, m_1^+$ and $m=\sum_{i=1}^{\ell}m_i^++~2, \dots, \sum_{i=1}^{\ell+1} m_{i}^+$ when $\ell=1, \ldots, {j-1}$. 
\end{itemize}
Let $G(s, t)$ be given by the expression \eqref{functionG}. Observe that for any $j \in \mathbb{N}$ the tangent to the curve $\mathcal{I}_p(t)$ at $t_{\sum_{\ell=1}^{j}m_\ell^+}=1$ is the vertical line $t=1$ which meets the tangent line at $t_{\sum_{\ell=1}^{j}m_{\ell}^+-1}$ at the point $\left(1, G(t_{\sum_{\ell=1}^{j}m_{\ell}^+-1}, 1)\right)$. Moreover, 
\[
\mathcal{I}_p(t)\geq \gamma_{t_{\sum_{\ell=1}^{j}m_{\ell}^+-1}}(t) \qquad \text{for} \qquad t \in [t_{\sum_{\ell=1}^{j}m_{\ell}^+-1}, 1]
\]
where $\gamma_s$ is given by \eqref{tangent}.
Hence,
\begin{align*}
 \int_{t_{\sum_{\ell=1}^{j}m_{\ell}^+-1}}^1 \cos\left(\frac{k\pi}{2}\mathcal{I}_p(t)\right) \mathrm{d}t >\mathcal{J}_{k, p}^{(1)}(t_{\sum_{\ell=1}^{j}m_{\ell}^+-1}, 1).
\end{align*}
and the proof is complete.
\end{proof}

\begin{remark}
If $k=3$ (so $j=1$) the formula above collapse to the simple expression
\begin{align*}
    \widehat{\mathrm{s}_p}(3)>\frac{4}{3\pi}&\left[\sum_{m=1}^{m_1^-}\mathcal{I}_{3,p}(x_{m-1},x_m)+\sum_{m=2}^{m_1^+-1}\mathcal{J}^{(1)}_{3,p}(t_{m-1},t_m)+ \right.\\
& \left. \sum_{m=2}^{m_1^+-1}\mathcal{J}^{(2)}_{3,p}(t_{m-1},t_m)+\mathcal{J}_{3,p}^{(1)}(t_{m_1^+-1},1)\right].
\end{align*}
The following table shows numerical lower bounds for $\widehat{\mathrm{s}_p}(3)$ whenever $p \in (1, \lambda]$.
\end{remark}
\centerline{
\begin{tabular}{c|c|c|c}
$\lambda$ & $m_1^-$ & $m_1^+$  & $\widehat{\mathrm{s}_p}(3)$ lower bound\\  \hline
1.5 & 2 & 3 & 0.0692320\\
1.5 & 3 & 3 & 0.0912921\\
1.5 & 4 & 3  & 0.0996541 \\
1.9 & 3 & 3 & 0.00534857\\
\end{tabular}}


\subsection{Towards analytic estimates for $\widehat{\mathrm{s}_p}(k)$ when $k\equiv_4 1$} \label{forckong1}
Let $k=4j-3$ for $j \in \mathbb{N}$. The function $\cos\big(\frac{k\pi}{2}\mathcal{I}_p(u)\big)$ with $u\in[0, 1]$ is monotonically decreasing in $j$ disjoint segments \[[\tilde{y}_i,\tilde{x}_i]  \qquad i=1,\ldots,j\] and it is monotonically increasing in $j-1$ disjoint segments \[[\tilde{s}_i,\tilde{t}_i] \qquad i=1,\ldots,j-1,\] so that
\[
      [0,1]=\left(\bigcup_{i=1}^{j} [\tilde{y}_i,\tilde{x}_i]  \right)
\cup
\left(\bigcup_{i=1}^{j-1} [\tilde{s}_i,\tilde{t}_i]  \right)
\]
where $\tilde{y}_1=0$, $\tilde{x}_{j}=1$, $\tilde{s}_i=\tilde{x}_i$ and $\tilde{y}_{i+1}=\tilde{t}_{i}$. The minimum turning points are such that 
\[
      \mathcal{I}_p(\tilde{x}_i)=\frac{4m-2}{k}\qquad \text{for }m=1,\ldots,j-1
\]
and the maximum turning points are such that
\[
    \mathcal{I}_p(\tilde{t}_i)=\frac{4m}{k}\qquad \text{for }m=1,\ldots,j-1.
\]

We partition each one of these segments into sets of quadrature points as follows. Let $\{m_{i}^-\}_{i=1}^{j},\, \{m_{i}^+\}_{i=1}^{j-1}\subset \mathbb{N}\setminus\{1\}$. Set 
\begin{align*}
&x_0=\tilde{y}_1=0,\quad x_{m_1^-}=\tilde{x}_1, \quad  x_{1+m_{1}^-}=\tilde{y}_2,\\ 
&x_{\sum_{\ell=1}^i m_\ell^-} =\tilde{x}_i,  & (i=2,\ldots,j)\\
&x_{1+\sum_{\ell=1}^i m_\ell^-}=\tilde{y}_{i+1},  & (i=2,\ldots,j-1)\\
& x_{\sum_{\ell=1}^{j} m_\ell^-}=\tilde{x}_{j}=1, \\
&t_1=\tilde{s}_1,\quad t_{m_1^+}=\tilde{t}_1,\quad t_{1+m_{1}^+}=\tilde{s}_2,\\
& t_{\sum_{\ell=1}^i m_\ell^+} =\tilde{t}_i, & (i=2,\ldots,j-1)\\ 
&t_{1+\sum_{\ell=1}^i m_\ell^+}=\tilde{s}_{i+1}, & (i=2,\ldots,j-2)
\end{align*}
We consider an increasing sequence of quadrature points 
\begin{align*}
&0\leq \cdots<x_{m-1}<x_{m}<\cdots\leq 1 &(m=1,\ldots,\sum_{\ell=1}^{j} m_\ell^-) \\
&0< \cdots<t_{m-1}<t_{m}<\cdots<1 &(m=2,\ldots,\sum_{\ell=1}^{j-1} m_\ell^+)
\end{align*}
such that
\begin{gather*}
   \left\{x_{1+\sum_{\ell=1}^{i-1} m_\ell^-}<\cdots<
   x_{\sum_{\ell=1}^{i} m_\ell^-}   \right\} \subset
[\tilde{y}_{i},\tilde{x}_i] \qquad \text{and} \\
\left\{t_{1+\sum_{\ell=1}^{i-1} m_\ell^+}<\cdots<
   t_{\sum_{\ell=1}^{i} m_\ell^+}    \right\} \subset
[\tilde{s}_{i},\tilde{t}_i].
\end{gather*}

\begin{mylem}\label{a19lowerb}
Let $p>1$ and $k=4j-3$ where $j\in \mathbb{N}$. For $k>1$
\begin{align*}
 \widehat{\mathrm{s}_p}(k)>\frac{4}{k\pi} &\left[ \sum_{m=1}^{m_1^-} \mathcal{I}_{k,p}(x_{m-1}, x_m)+\sum_{\ell=1}^{j-1} \sum_{m=\sum_{i=1}^{\ell}m_i^-+2}^{\sum_{i=1}^{\ell+1} m_{i}^-} \mathcal{I}_{k, p}(x_{m-1}, x_m)+ \right. \\ & \left.+\sum_{m=2}^{m_1^+} \Big(\mathcal{J}_{k, p}^{(1)}(t_{m-1}, t_m)+ \mathcal{J}_{k, p}^{(2)}(t_{m-1}, t_m)\Big)\right. \\ & \left. +\sum_{\ell=1}^{j-2}\sum_{m=\sum_{i=1}^{\ell}m_i^++2}^{\sum_{i=1}^{\ell+1}m_{i}^+}
\Big(\mathcal{J}_{k, p}^{(1)}(t_{m-1}, t_m)+\mathcal{J}_{k, p}^{(2)}(t_{m-1}, t_m)\Big) \right].
\end{align*}
\end{mylem}
\begin{proof}
The proof is similar to that of Lemma~\ref{a13lowerb}.
\end{proof}

Numerically we have the following lower bounds for $\widehat{\mathrm{s}_p}(9)$ whenever $p \in (1, \lambda]$.

\centerline{
\begin{tabular}{c|c|c|c|c|c|c}
$\lambda$ & $m_1^-$ & $m_1^+$ & $m_2^-$ &$ m_2^+$ & $m_3^-$ & $\widehat{\mathrm{s}_p}(9)$ lower bound\\  \hline
1.5 & 4 & 5 & 5 & 4 & 2 &  $8.76881\times 10^{-6}$\\
1.5 & 5 & 5 & 5 & 4 & 2 &  $8.35771\times 10^{-5}$
\end{tabular}}


\section*{Acknowledgements}
This research initiated during a visit of LB to The Ohio State University in December~2016. He is kindly grateful to Jan Lang and Boris Mityagin for their helpful comments and the financial support provided.


\begin{thebibliography}{10}

\bibitem{BBCDG2006}
{\sc P.~Binding, L.~Boulton, J.~{\v C}epi{\v c}ka, P.~Dr{\'a}bek, and P.~Girg},
  {\em Basis properties of eigenfunctions of the {$p$}-{L}aplacian}, Proc.
  Amer. Math. Soc., 134 (2006), pp.~3487--3494.

\bibitem{BL2014}
{\sc L.~Boulton and G.~J. Lord}, {\em Basis properties of the $p,q$-sine
  functions}, Proc. R. Soc. A, 471 (2015).

\bibitem{BushellEdmunds2012}
{\sc P.~J. Bushell and D.~E. Edmunds}, {\em Remarks on generalized
  trigonometric functions}, Rocky Mountain J. Math., 42 (2012), pp.~25--57.

\bibitem{EdmundsGurkaLang2012}
{\sc D.~E. Edmunds, P.~Gurka, and J.~Lang}, {\em Properties of generalized
  trigonometric functions}, J. Approx. Theory, 164 (2012), pp.~47--56.

\bibitem{Grad2007}
{\sc I.~S. Gradshteyn and I.~M. Ryzhik}, {\em Table of Integrals, Series and
  Products}, Elsevier/Academic Press, Amsterdam, seventh~ed., 2007.

\bibitem{HLS1997}
{\sc H.~Hedenmalm, P.~Lindqvist, and K.~Seip}, {\em A {H}ilbert space of
  {D}irichlet series and systems of dilated functions in ${L^2(0,1)}$}, Duke
  Mathematical Journal, 86 (1997), pp.~1--37.

\bibitem{HLS1999}
\leavevmode\vrule height 2pt depth -1.6pt width 23pt, {\em Addendum to: A
  {H}ilbert space of {D}irichlet series and systems of dilated functions in
  ${L^2(0,1)}$}, Duke Mathematical Journal, 99 (1999), pp.~175--178.

\bibitem{LangEdmunds2011}
{\sc J.~Lang and D.~Edmunds}, {\em Eigenvalues, embeddings and generalised
  trigonometric functions}, vol.~2016 of Lecture Notes in Mathematics,
  Springer, Heidelberg, 2011.

\bibitem{Mit2017}
{\sc B.~Mityagin}, {\em Systems of dilated functions: completeness, minimality,
  basisness}, Funct. Anal. Appl., 51 (2017), pp.~236--239.

\bibitem{Singer1970}
{\sc I.~Singer}, {\em Bases in {B}anach spaces. {I}}, Springer-Verlag, New
  York-Berlin, 1970.
\newblock Die Grundlehren der mathematischen Wissenschaften, Band 154.

\end{thebibliography}
\end{document}